\documentclass[preprint,12pt]{elsarticle}

\usepackage{amssymb}
\usepackage{amsthm}
\usepackage{subcaption}

\theoremstyle{definition}

\theoremstyle{remark}

\def \bc{\begin{center}}
\def \ec{\end{center}}

\begin{document}

\begin{frontmatter}

\title{Subspace method based on neural networks  for solving the partial differential equation in weak form}

\author[mymainaddress,mythirdaddress]{Pengyuan Liu}\ead{liupengyuan23@gscaep.ac.cn}
\author[mymainaddress,mythirdaddress]{Zhaodong Xu}\ead{ xuzhaodong\_math@163.com}
\author[mymainaddress,mysecondaryaddress]{Zhiqiang Sheng  \corref{cor1}}\ead{sheng\_zhiqiang@iapcm.ac.cn}

\cortext[cor1]{Corresponding author.}

\address[mymainaddress]{National Key Laboratory of Computational Physics, Institute of Applied Physics and Computational Mathematics,  Beijing, 100088, China}
\address[mysecondaryaddress]{HEDPS, Center for Applied Physics and Technology, and College of Engineering, Peking University, Beijing, 100871, China}
\address[mythirdaddress]{Graduate School of China Academy of Engineering Physics, Beijing 100088, China}

\begin{abstract}
We  present a subspace method based on neural networks   for solving the partial differential equation in weak form with high accuracy.
The basic idea of our method is to use some functions based on neural networks  as base functions to span a subspace, then find an approximate solution in this subspace.
Training base functions and finding an approximate solution can be separated, that is
different  methods can be  used to train these base functions, and different methods can also be used to find an approximate solution.
In this paper, we find an approximate solution of the partial differential equation in the weak form.
Our method can achieve high accuracy with low cost of training.
Numerical examples show that the cost of training these base functions
is low, and only one hundred to two thousand epochs are needed for most  tests.
The error of our method can  fall below the level of $10^{-7}$ for some tests.
The proposed method has the better performance in terms of the accuracy and computational cost.

\end{abstract}


\end{frontmatter}

\date{}



\section{Introduction}\label{introduce}

Due to the rapid development of machine learning methods,
the  method based on neural networks attracts more and more attention.
Many  numerical methods based on neural networks have been proposed for solving the partial differential equation.
Some of these methods are based on the strong form of the partial differential equation,
such as physical information neural networks (PINN)\cite{Raissi-19-jcp}, deep Galerkin method (DGM)\cite{Sirignano-18-jcp}, methods based on extreme learning machine (ELM)\cite{Dong-21-cmame} and random feature methods (RFM) \cite{Chen-22-jml}. And others are based on the weak form of the partial differential equation, such as deep Ritz method (DRM)\cite{E-18}, and weak adversarial networks (WAN)\cite{Zang-20-jcp}.

The loss functions of PINN  are based on  the mean squared error  consisting of information about the partial differential equation in  strong form as well as the initial boundary conditions on some collocations points \cite{Jagtap-20-cmame,Jagtap-22-jcp,Liu-24-jsc,Patel-22-jcp,Yao-23-amame,Ying-23-amm}.
The loss functions of DGM are based on the $L^2$-norm
of the residuals of partial differential equation in  strong form.
For these ELM-based methods   and random feature methods,  these parameters of neural networks are generated randomly and do not need to be updated,  the approximate solution is given by solving a least-squares problem \cite{Chi-24-cmame, Li-03087,Lyu-22-jcp,Shang-23,Sun-24-jcam}.
The loss function of the deep Ritz method is  based on the energy functional corresponding to the weak form of partial differential equation.
The weak adversarial networks method constructs a loss function by minimizing an operator norm
induced from the weak form of partial differential equation.

The deep Ritz Method \cite{E-18} uses the variational form of the partial differential equation as the loss function, and minimizes the loss function to update these parameters of neural networks.
The deep Nitsche Method is proposed to deal with the essential boundary conditions in \cite{Liao-21-cicp}.
This method imposes the boundary conditions in a nonconforming way as the penalty method.
A penalty-free neural network method for  solving a class of second-order boundary-value problems on complex geometries is proposed in \cite{sheng-21-jcp}. In which, the original problem is reformulated to a weak form so that the evaluations of high-order derivatives are avoided.
An approach   based on randomized neural networks and the Petrov-Galerkin method is developed in \cite{Shang-23}, in which finite element basis functions can be used as the test functions.
A local randomized neural networks method with discontinuous Galerkin  methods  is  proposed in \cite{Sun-24-jcam}, and the discontinuous Galerkin  method is used to glue the solutions on sub-domains together.
Some methods based on tensor neural networks for  solving high-dimensional partial differential equation are proposed in \cite{Wang-doctor,Wang02754, Wang02732}.

However, the accuracy of  most  methods is unsatisfactory, and the training  cost  is extremely high.
It is a challenge   to compete with traditional methods for low dimensional problems.
A subspace method based on neural networks (SNN) for solving the partial differential equation in the strong form is proposed in \cite{xu-24}.
The basic idea  is to use some functions based on neural networks  as base functions to span a subspace, then find an approximate solution in this subspace.
Two special algorithms in  the strong form of partial differential equation are designed.
This method  achieves high accuracy with low cost of training.

In this paper, we present a subspace method based on neural networks  for solving the partial differential equation in the weak form (SNNW).
First, we train these base functions in the subspace layer. These parameters of neural networks are updated by minimizing the loss function. The loss function can be based on different form, such as the strong form and the weak form. This means that these parameters can be updated by using various methods.
We do not restrict the use of specific methods for training these base functions, and can use PINN, DGM or DRM to train these base functions in this step.
After obtaining these base functions, we  find an approximate solution of the partial differential equation in the weak form.
Similar to \cite{xu-24}, the loss function can include only the information of PDE and do not include the information of the initial boundary condition.

Numerical examples show that the cost of training these base functions
is low, and only one hundred to two thousand epochs are needed for most tests.
The error of our method can  fall below the level of $10^{-7}$ for some tests.
The performance of our method significantly surpasses the performance of existing methods based on the weak form of PDE in terms of the accuracy and computational cost.

The remainder of this paper is organized as follows.
In section 2, we   describe  the subspace method based on neural networks for solving the partial differential equation in the weak form.
In section 3, we present some numerical examples  to test the performance of our method.   At last, we give some conclusions.

\section{SNN for solving the PDE in weak form} 
Consider the  following equation:
 \begin{eqnarray}
\mathcal{A}u(\bold{x})&=&f(\bold{x})  \qquad \ \ \mbox{in}\ \ \Omega,
\label{deq1} \\
\mathcal{B} u(\bold{x}) &=&g(\bold{x})  \qquad \ \  \mbox{on} \ \
\partial \Omega,  \label{deq2}
\end{eqnarray}
where $\bold{x}=(x_1,x_2,\cdots, x_d)^T$,   $\Omega$ is a bounded domain in ${\textrm{R}}^d$, $\partial \Omega$ is the boundary of $\Omega$, $\mathcal{A}$ and $\mathcal{B}$ are the differential operators, $f$ and $g$ are given functions.

\subsection{Neural networks architecture}
In this paper, we use the same neural networks architecture as \cite{xu-24}. For completeness, we briefly describe the neural networks architecture.

\begin{figure}[htbp]
\begin{center}\includegraphics[width = 10cm]{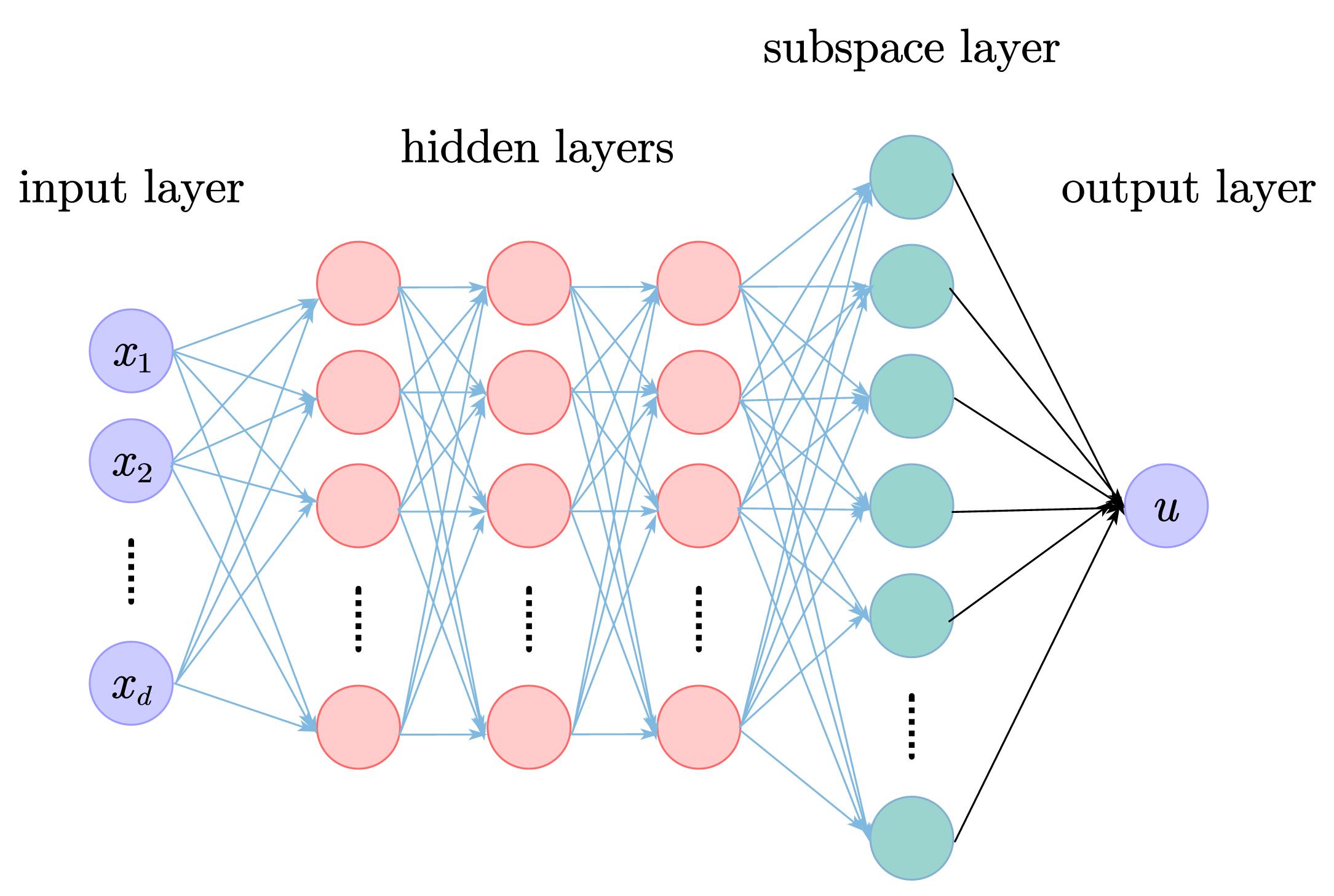}
\caption{The neural networks architecture.} \label{network}
\end{center}
\end{figure}

The neural networks architecture consists of four layers, which include input layer, hidden layers, subspace layer and  output layer.    Figure \ref{network} illustrates the neural networks architecture.

Let $M$ be the dimension of subspace in the subspace layer, and $\varphi_j$($j=1,2,\cdots,M$) be base functions of subspace,
 and $\omega_j$($j=1,2,\cdots,M$) be some coefficients related to base functions.
 Denote $\varphi = (\varphi_1,\varphi_2,\cdots, \varphi_M)^T$ and $\omega =(\omega_1, \omega_2, \cdots, \omega_M)^T$.
Let
$\theta$ be the set of  parameters in neural networks, $u(x;\theta, \omega)$ is the output with respect to input $x$ with parameters $\theta$ and $\omega$.

\subsection{The training of base functions}\label{basefunction}

We train the base functions of subspace such that the subspace has effective approximate capability to the solution space of equation.

In this step, we can use different methods to train the base functions of subspace.
We describe how to train the base functions by using PINN, DGM and DRM.

\subsubsection{The base functions based on PINN} \label{PINN}
To avoid introducing this penalty parameter,  we  define the following loss function which  contains only the PDE loss term.
\begin{equation}
\mathcal{L_{D}}=\frac1N \sum_{j=1}^N \left(\mathcal{A} u(\bold{x}_j;\theta,\omega) -f(\bold{x}_j)\right)^2.
\label{loss-d}
\end{equation}
It is obvious  that this loss function  includes only the information of PDE itself, and omit the information of initial boundary conditions.

First, we  train the parameter $\theta$ and $\omega$ by  minimizing the loss function $\mathcal{L}_{D}(\bold{x};\theta,\omega)$ in (\ref{loss-d}).
In order to balance the accuracy and efficiency, it is not necessary to solve
 the minimization problem accurately. Let ${\mathcal{L}_{D}}_0(\bold{x};\theta,\omega)$ be the initial loss.
 The training process stops if the following condition is satisfied,
\begin{equation}
\frac{\mathcal{L}_{D}(\bold{x};\theta,\omega)}{{\mathcal{L}_{D}}_0(\bold{x};\theta,\omega)} \le \varepsilon.
\label{loss-dstop}
\end{equation}
This implies that  the training process stops if the  loss decreases to a certain level.
To prevent excessive  training epochs, we introduce the maximum number of training epochs $N_{max}$.
If the number of training epochs reaches  $N_{max}$,  the training process also stops.

\subsubsection{The base functions based on DGM}\label{DGM}

Define  the loss function as follows:
\begin{equation}
\mathcal{L}_I=\| \mathcal{A} u(\bold{x};\theta,\omega) -f(\bold{x})\|_{L^2(\Omega)}^2.
\label{loss-i}
\end{equation}
We can see that the loss function includes  only the information of PDE, and omit the information of initial boundary conditions.

At first, we   train the parameter $\theta$ and $\omega$ by  minimizing the loss function $\mathcal{L}_{I}(\bold{x};\theta,\omega)$ in (\ref{loss-i}).
Let ${\mathcal{L}_{I}}_0(\bold{x};\theta,\omega)$ be the initial loss.
 The training process stops if the following condition is satisfied,
\begin{equation}
\frac{\mathcal{L}_{I}(\bold{x};\theta,\omega)}{{\mathcal{L}_{I}}_0(\bold{x};\theta,\omega)} \le \varepsilon.
\label{loss-istop}
\end{equation}
This implies that  the training process stops if the  loss decreases to a certain level.
If the epochs reach $N_{max}$,  the training process also stops.

\subsubsection{The base functions based on DRM}\label{DRM}
Define  the loss function as follows:
\begin{equation}
\mathcal{L}_{R}=\frac12 a( u(\bold{x};\theta,\omega),u(\bold{x};\theta,\omega)) -(f,u(\bold{x};\theta,\omega)),
\label{loss-drm}
\end{equation}
where $a(\cdot,\cdot)$ is a bilinear form that depends on   the  operators $\mathcal{A}$ and $\mathcal{B}$. We will describe this bilinear form in the next subsection.
Different from the loss function of PINN and DGM, the minimum value of the loss function $\mathcal{L}_{R}$ is not 0. Hence, we can not use the relative loss to stop the training.
In this paper, we can stop the training after the training epochs reach a given number.

\subsection{The weak form}

The equation (\ref{deq1}) with the boundary condition (\ref{deq2}) has the weak form: Find $u \in U$ such that
\begin{eqnarray}
a(u,v) = (f,v), \quad \forall v \in V, \label{weakform}
\end{eqnarray}
where $U$ and $V$ are two Hilbert spaces, $a(\cdot,\cdot)$ is a bilinear form that depends on   the  operators $\mathcal{A}$ and $\mathcal{B}$.

We take the Poisson equation as the model problem to give the weak form.
Consider the   Poisson equation  on the domain $\Omega=(0,1) \times(0,1)$,
\begin{eqnarray}
        -\Delta u&=&f,  \quad \mbox{in} \ \Omega, \label{poisson}\\
    u&=&0,  \quad \mbox{on} \ \partial \Omega. \label{poisson-boundary}
\end{eqnarray}
From (\ref{poisson}), there is
\begin{eqnarray}
-\int_{\Omega} \Delta u v dx&=\int_{\Omega} f v dx, \,\,\,\,\forall v\in H_0^{1}(\Omega).
\end{eqnarray}
By the Green's formula, we have
\begin{eqnarray}
-\int_{\partial\Omega}  \frac{\partial u}{\partial\mathbf{n}}  v dx+\int_{\Omega}  \nabla u\cdot\nabla v dx =\int_{\Omega} f v dx, \,\,\,\,\forall v\in H_0^{1}(\Omega).
\end{eqnarray}
It follows
\begin{eqnarray}
\int_{\Omega}  \nabla u\cdot\nabla v dx=\int_{\Omega} f v dx, \,\,\,\,\forall v\in H_0^{1}(\Omega).
\end{eqnarray}
Hence, we can get the  following  weak form:
\begin{eqnarray}
a(u,v)=(f,v), \,\,\,\,\forall v\in H_0^{1}(\Omega), \label{weak-poisson}
\end{eqnarray}
where
\[
a(u,v)=\int_{\Omega}  \nabla u\cdot\nabla v dx,
\]
and
\[
(f,v)=\int_{\Omega} f v dx.
\]

In order to avoid the penalty on boundary condition, we simply use the following method to treat the Dirichlet boundary condition.
For the problem (\ref{poisson})-(\ref{poisson-boundary}), we multiply each $\bar{\varphi}_{j}$ by the limit function $h:x(1-x)y(1-y)$ to obtain the homogeneous trial function subspace as follows:
\[
{\varphi}_{j} =x(1-x)y(1-y)\bar{\varphi}_{j},\quad j=1,\cdots, M.
\]
where $\varphi_{j}$ and $\bar{\varphi}_{j}$ are obtained by the training process of subspace layer in section \ref{basefunction}.
It should be pointed out that there are many choices of the
limiting function $h$ as long as they satisfy the condition of being zero on $\partial \Omega$ and nonzero in $\Omega$.

Define the finite dimensional space  $V_h=span\{{\varphi}_{1}, {\varphi}_{2}, \cdots, {\varphi}_{M} \}$.
The corresponding weak form: 
Find $u_h \in V_h$ such that
\begin{eqnarray}\label{weakform}
a(u_h,v_h)=(f,v_h), \,\,\,\,\forall v_h\in V_h.  \label{weak-helm}
\end{eqnarray}
Note that $u_h$ can be expressed as
\begin{eqnarray}
u_{h}=\sum^{M}_{j=1}\omega_j {\varphi}_{j}. \label{uh}
\end{eqnarray}
Inserting (\ref{uh}) into (\ref{weak-helm}), and setting $v_h = {\varphi}_i, i=1, \cdots,M$,  we obtain
\begin{eqnarray}
\sum^{M}_{j=1}\omega_{j}a({\varphi}_{j},{\varphi}_{i})=(f,{\varphi}_{i}),\quad i=1,2,\cdots, M.\label{weak-equation}
\end{eqnarray}
From (\ref{weak-equation}), we can get a linear system
\begin{eqnarray}
A\omega=b, \label{linear-system}
\end{eqnarray}
where $A$ is a $M\times M$ matrix with $A_{ij} = a({\varphi}_{j},{\varphi}_i)$ and $b$ is a $M\times 1$ vector with $b_i = (f,{\varphi}_i)$.
Solve the linear system (\ref{linear-system}) to obtian $\omega$, and then the approximate solution is obtained by (\ref{uh}).

We summarize the main steps of our method as follows:

\begin{tabular}{l}\hline
 {\bf Algorithm:}  SNNW\\
\hline
1. Initialize nerual networks architecture.\\

2. Generate randomly $\theta$ and  $\omega$.\\

3. Update parameter $\theta$ and  $\omega$ by minimizing the loss function $\mathcal{L}(\bold{x};\theta,\omega)$  \\
   \ \quad  until stop criteria is satisfied. \\

4. Obtain the base functions of the subspace $\varphi_1$, $\varphi_2$, $\cdots$, $\varphi_M$. \\

5. Solve the algebraic system (\ref{linear-system}) resulted from (\ref{weakform}) in weak form to \\
   \ \quad update $\omega$.\\

6. Obtain an approximate solution $u_h$.  \\
\hline

\end{tabular}

\vspace{0.3cm}

Different loss function leads different method.
SNNW based on the base functions in Section \ref{PINN} is denoted as SNNW-P,
SNNW based on the base functions in Section  \ref{DGM} is denoted as SNNW-G,
and SNNW based on the base functions in Section \ref{DRM} is denoted as SNNW-R.

\section{Numerical results}
\label{Numerical results}
In this section, we use some numerical experiments to test our method. In Section \ref {Helmholtz equation}, we test the performance of SNNW by solving the one-dimensional Helmholtz equation.  In Section \ref {Poisson equation}, we test our method for solving the  Poisson equation. In Section \ref{Anisotropic diffusion equation}, we  solve the  diffusion equation with strongly anisotropic ratios.

We employ the  framework PyTorch \cite{paszke-19-ansps}, and  all variable data types are set to float64.
We maintain consistent settings for all numerical tests in this paper. Specifically, we use a feedforward fully connected neural network with four hidden layers, each containing 100 neurons. The activation function is the Tanh function,  the subspace dimension is  300, the optimizer is Adam and the settings of optimizer  are kept at their default values. Neural network parameters are randomly generated by using the Xavier method.  In order to accurately reproduce numerical results, the seed for generating random numbers  in all numerical experiments is set to 1.

For SNNW-P and SNNW-G, the training process stops  when the relative  loss is less than $\varepsilon$ or  the epochs reach   $N_{max}$. For SNNW-R, the training process stops  when the epochs reach  2000.
In all tests, we take $\varepsilon=1e-3$ and $N_{max}=5000$.
 
We compare our method with the existing methods, including PINN, DGM and DRM.
For comparison purposes, identical network architectures, parameter settings, and initialization methods are employed for PINN, DGM and DRM, with the Adam optimizer training through 50000 epochs.

We  use the following relative $L^2$ error to evaluate the accuracy,
\begin{eqnarray*}
\|e\|_{ {L}^2}=\frac{\sqrt{\sum_{i=1}^N\left|u_h\left(X_i\right)
-u^*\left(X_i\right)\right|^2}}{\sqrt{\sum_{i=1}^N\left|u^*\left(X_i\right)\right|^2}},
\end{eqnarray*}
where $ u^*$ is the exact solution.

\subsection{Helmholtz equation} \label{Helmholtz equation}
Consider the 1D Helmholtz equation with the Dirichlet boundary conditions as follows:
\begin{eqnarray*}
&-u_{xx}+\lambda u=f,\,\,a<x<b,\\
&u(a)=0,\ u(b)=0.
\end{eqnarray*}
The exact solution  is taken as
\begin{eqnarray*}\label{hzsolution}
u(x)=\sin{(3\pi x)}+\cos{(4\pi x +\frac{1}{2}\pi)},
\end{eqnarray*}
where $a=0\,,b=2,\,\lambda=1$.

To obtain the integration in weak form of SNNW directly, the complex Guassian points are chosen in training process. The complex Gaussian
quadrature formula is applied to segment $[a, b]$ into 100 sub-intervals, with
each sub-interval containing 10 points.

Table \ref{helmholtz error different method} shows the relative $L^2$ error for different methods.
Obviously, the errors of PINN, DGM and DRM  are between $10^{-3}$ and $10^{-4}$ after 50000 training epochs.
We can see  that SNNW-P achieves relative $L^2$ error of 3.50e-06 with only 250 epochs, SNNW-G achieves relative $L^2$ error of 6.73e-06 with only 273 epochs, and SNNW-R achieves relative $L^2$ error of 5.65e-06 with 2000 epochs. This demonstrates that the performance of SNNW has the superiority in the accuracy and training cost.

\begin{table}[!htbp]
\caption{The errors and epochs of different methods for the Helmholtz equation}\label{helmholtz error different method}
\begin{center}
\small
\begin{tabular}{cccc}\hline
Method & $\|e\|_{L^2}$ & epochs\\
\hline
PINN &  5.63e-03 & 50000 \\
DGM  &  2.29e-04 & 50000 \\
DRM  &  4.51e-03 &  50000\\
SNNW-P & 3.50e-06 & 250 \\
SNNW-G & 6.73e-06 & 273 \\
SNNW-R & 5.65e-06  & 2000\\
\hline
\end{tabular}
\end{center}
\end{table}

\begin{figure}[htbp]
\begin{center}\includegraphics[width = 6cm]{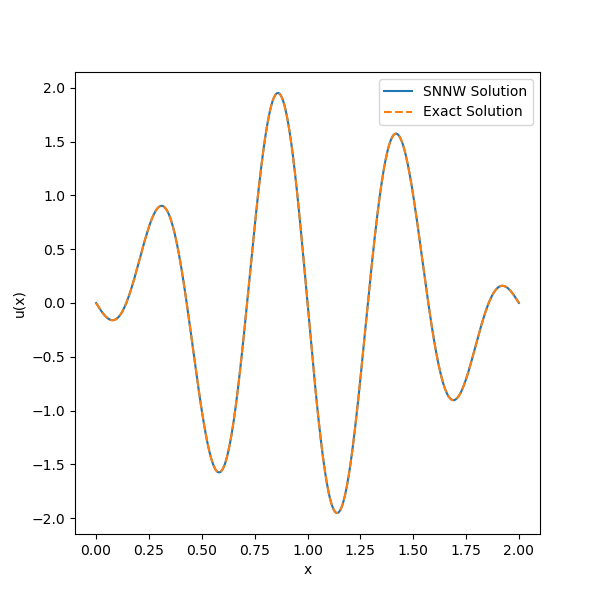}
 \caption{Solution obtained by SNNW for Helmholtz Equation.} \label{hz_exact_vs_appro}
\end{center}
\end{figure}

\begin{figure}[htbp]
\centering
\subcaptionbox{SNNW-P}{\includegraphics[width=0.3\linewidth]{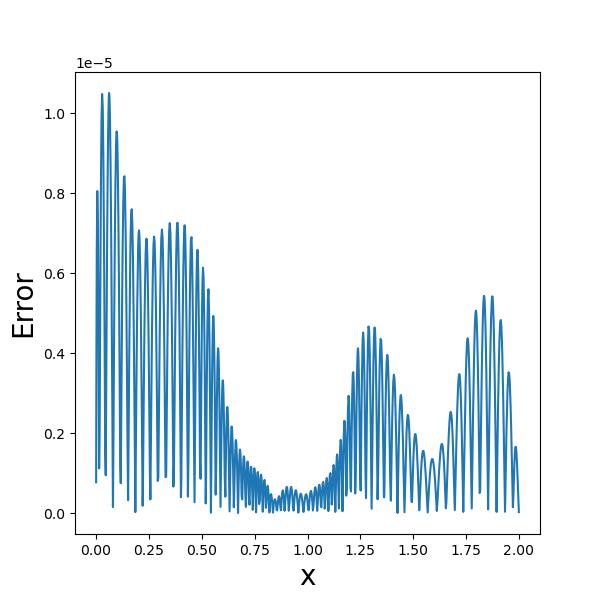}}
\hspace{0.3cm} 
\subcaptionbox{SNNW-G}{\includegraphics[width=0.3\linewidth]{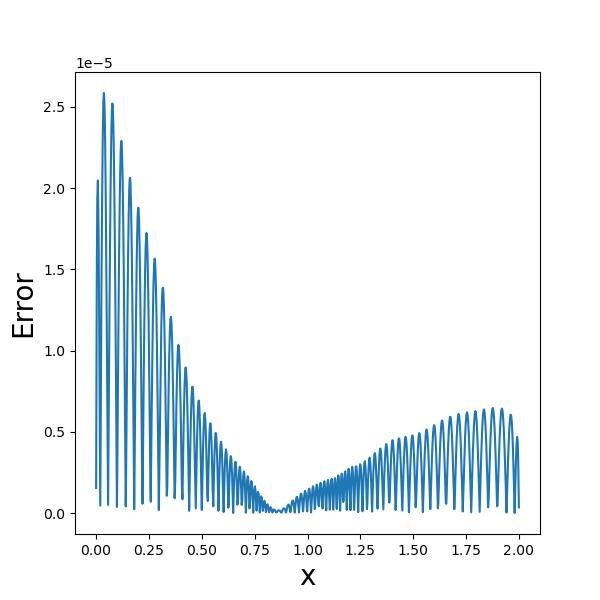}}
\hspace{0.3cm}  
\subcaptionbox{SNNW-R}{\includegraphics[width=0.3\linewidth]{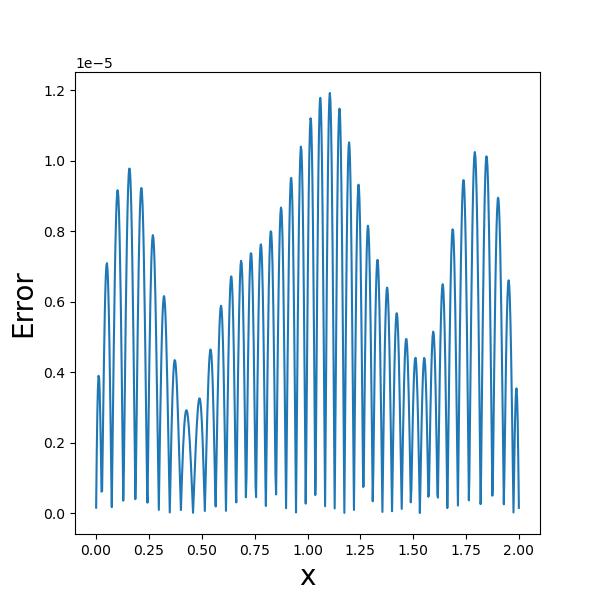}}
\caption{Point-wise errors of  SNNW-P, SNNW-G and SNNW-R for Helmholtz equation.}
\label{ hz_Point-wise_error}
\end{figure}

The  numerical solutions obtained by our SNNW method are shown in Figure \ref{hz_exact_vs_appro}. Figure \ref{ hz_Point-wise_error} illustrates the point-wise errors of our methods.

The numerical results of SNNW-P and SNNW-G are displayed for various numbers of sampling points and subspace dimension in Tables \ref{hzSNNWP} and   \ref{hzSNNWG} .
The network setup is uniform with four hidden layers, each with 100 neurons and adopts the composite Gaussian quadrature rule, sampling 10 points per sub-interval and incrementing sub-interval
counts for accuracy enhancement.
With the sampling point number of 60 and 80, the accuracy is not satisfied.  This means that there needs more sampling points to obtain better accuracy. Increasing the number of sampling point  to 200, SNNW-P and SNNW-G can achieve rapidly the accuracy of $10^{-5}$.
 As the number of sampling points and subspace dimension increases, both SNNW-P and SNNW-G can obtain the stable accuracy of $10^{-6}$. Moreover, the number of training epochs is between 200 and 3000.

\begin{table}[!htbp]
\caption{The errors and epochs of SNNW-P across various numbers of sampling points and subspace dimension $M$ for the Helmholtz equation.}\label{helmholtz error SNNWP points and subspace dimension }\label{hzSNNWP}
\begin{center}
\resizebox{\columnwidth}{!}{%
\small
\begin{tabular}{ccccccccc}
\hline
Points & $M$ & 60 & 80 & 100 & 300 & 500  \\
\hline
$60$ & $\|e\|_{ {L}^2}$ & 4.33e01  & 1.12e02  & 6.66e01 & 5.04e01  & 9.98e00   \\
 & epochs & 2329  & 399  & 1406  & 366   & 506    \\
$80$ & $\|e\|_{ {L}^2}$ & 3.85e00  & 1.45e01  & 4.15e01  & 8.71e-01  & 3.77e-03   \\
 & epochs & 994  & 721  & 1029  & 431   & 285    \\
 $100$ & $\|e\|_{ {L}^2}$ & 2.06e-03 & 4.69e-04 & 6.76e-03 & 8.85e00  & 2.79e-04   \\
 & epochs & 347 & 287 & 307 & 288  &  298  \\
 $200$ & $\|e\|_{ {L}^2}$ & 9.56e-05 & 4.46e-06 & 1.33e-05 & 1.13e-05 & 1.75e-06  \\
 & epochs & 275 & 261 &  574 & 309  &  280  \\
 $300$ & $\|e\|_{ {L}^2}$ &  9.40e-05 & 2.10e-05  & 2.78e-05  & 4.36e-06  & 1.09e-06   \\
 & epochs & 352 & 321  & 660  & 298   & 258    \\
 $500$ & $\|e\|_{ {L}^2}$ & 1.37e-04  &  4.64e-05 &  3.01e-05 & 5.94e-06 & 1.60e-06   \\
 & epochs & 267 & 430  & 807 &  276  & 249   \\
 $1000$ & $\|e\|_{ {L}^2}$ & 6.67e-05 & 7.94e-05 & 2.70e-05  & 3.40e-06  & 1.46e-06  \\
 & epochs & 295 & 420  & 676  & 282  & 250    \\
\hline
\end{tabular}
} 
\end{center}
\end{table}

\begin{table}[!htbp]
\caption{The errors and epochs of SNNW-G across various numbers of sampling points and subspace dimension $M$ for the Helmholtz equation.}\label{hzSNNWG}
\begin{center}
\resizebox{\columnwidth}{!}{%
\small
\begin{tabular}{ccccccccc}
\hline
Points & $M$ & 60 & 80 & 100 & 300 & 500  \\
\hline
$60$ & $\|e\|_{ {L}^2}$ & 7.85e01  & 6.25e01  & 6.31e01 & 4.94e01  & 9.43e-02   \\
 & epochs &  537 &  330  & 2202  &  396  &  354   \\
$80$ & $\|e\|_{ {L}^2}$ & 8.38e-01  & 2.52e01   &  4.66e01  &  1.91e01 & 3.13e-04   \\
 & epochs & 449  &  378 & 1211 &  583  & 260    \\
 $100$ & $\|e\|_{ {L}^2}$ & 1.84e-03 &  5.10e-01 & 8.88e-04 & 7.42e-04  & 1.48e-03  \\
 & epochs &  305 & 584 &  440 & 284  & 302  \\
 $200$ & $\|e\|_{ {L}^2}$ & 1.65e-04 & 7,67e-05 & 6.11e-05  & 2.38e-05  & 2.43e-06   \\
 & epochs &310 &  531&  425 &  245 &  300 \\
 $300$ & $\|e\|_{ {L}^2}$ & 8.62e-05 & 7.04e-05 & 1.86e-05   &  2.58e-06 & 3.05e-06   \\
 & epochs &  286  & 324  & 365  & 246  &  268   \\
 $500$ & $\|e\|_{ {L}^2}$ & 1.01e-04 & 3.11e-05   & 3.04e-05   &6.30e-06   & 3.32e-06  \\
 & epochs & 307 &  347 & 381  & 270 & 267   \\
 $1000$ & $\|e\|_{ {L}^2}$ &  1.23e-04 & 1.64e-05   & 1.92e-05   & 6.73e-06 & 3.19e-06    \\
 & epochs &  314  & 316  &  362 &  273  &  248  \\
\hline
\end{tabular}
} 
\end{center}
\end{table}
 Similar to SNNW-P and SNNW-G, SNNW-R has the same network setup and composite Gaussian quadrature rule. The numerical results of SNNW-R with given 2000 epochs are shown in  Table \ref{hzSNNWR} for various numbers of sampling points and subspace dimension.
 The accuracy is not satisfied when the number of sampling point number is small. Increasing the number of sampling points and subspace dimension, SNNW-R obtains the accuracy of $10^{-5}$.
 The number of training epochs of our methods is significantly less than that
required by PINN, DGM and DRM.

\begin{table}[!htbp]
\caption{The errors of SNNW-R with 2000 epochs across various numbers of sampling points and subspace dimension $M$ for the Helmholtz equation}\label{hzSNNWR}
\begin{center}
\resizebox{\columnwidth}{!}{%
\small
\begin{tabular}{ccccccccc}
\hline
Points & $M$ & 60 & 80 & 100 & 200 & 300 & 500\\
\hline
$60$ & $\|e\|_{L^2}$ &  3.00e01 & 3.59e01 & 2.81e01 & 6.51e00 & 6.30e-04 & 7.38e-04 \\
$80$ & $\|e\|_{L^2}$ & 1.38e-03 & 1.13e-02  & 7.52e-02  & 9.39e-03  &  7.58e-06  & 1.05e-05  \\
$100$ & $\|e\|_{L^2}$ & 1.50e-04 & 1.42e-04 & 1.18e-03  & 3.79e-05  & 6.10e-06 & 1.30e-05 \\
$200$ & $\|e\|_{L^2}$ & 9.06e-05 & 5.06e-06 & 5.00e-06  & 2.39e-06  & 5.64e-06  & 1.27e-05 \\
$300$ & $\|e\|_{L^2}$ &  8.98e-05 & 6.38e-06 & 5.04e-06  & 2.09e-06  & 5.74e-06 & 1.24e-05 \\
$500$ & $\|e\|_{L^2}$ & 1.02e-04 & 6.65e-06 &  5.06e-06 &  2.16e-06 & 5.53e-06 &  1.26e-05\\
$1000$ & $\|e\|_{L^2}$ & 1.00e-04 & 6.46e-06  & 5.07e-06  &  2.39e-06  & 5.65e-06   & 1.25e-05  \\
\hline
\end{tabular}
} 
\end{center}
\end{table}

\begin{figure}[htbp]
\centering
\subcaptionbox{}{\includegraphics[width=0.45\linewidth]{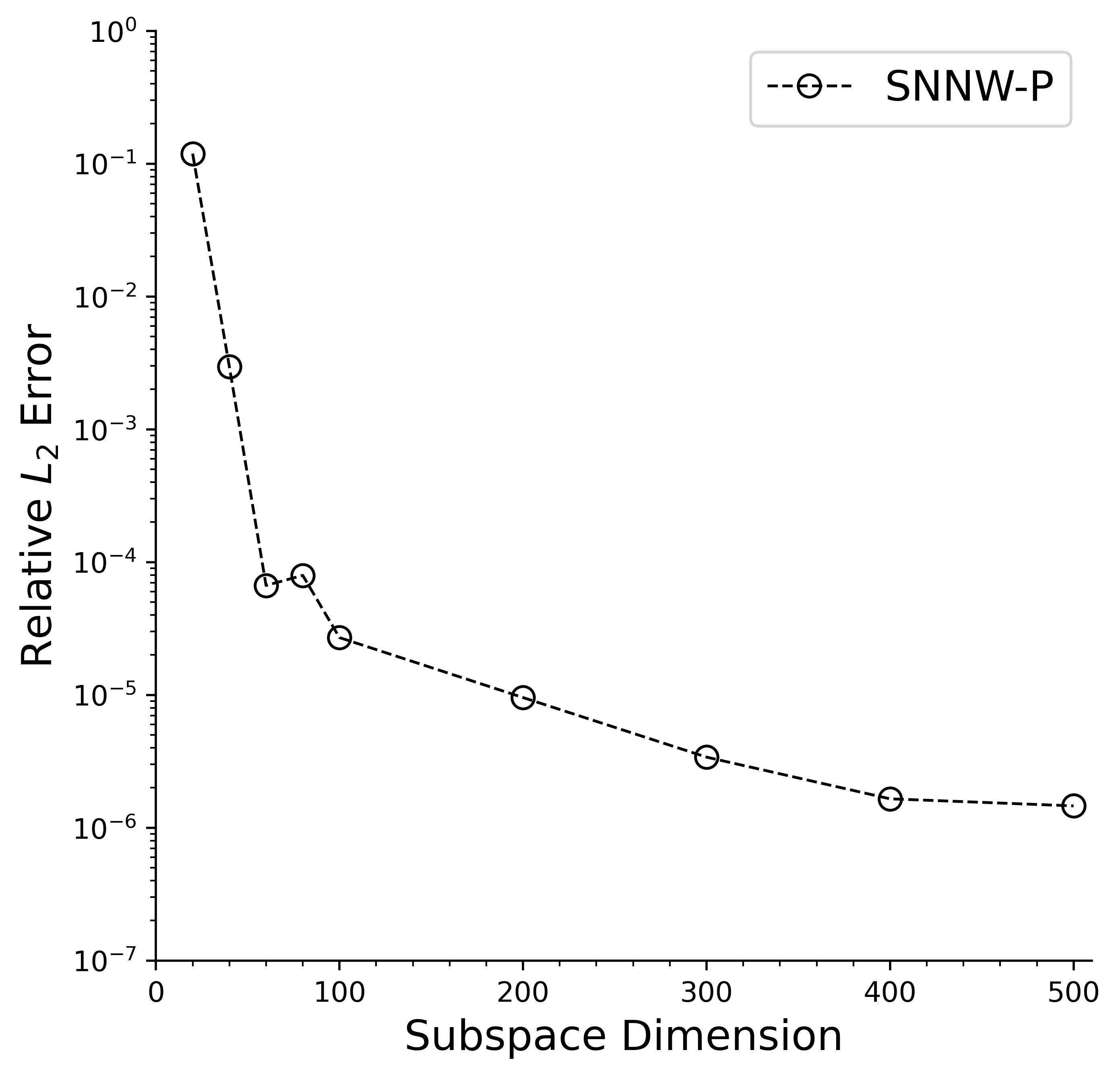}}
\hspace{0.3cm}  
\subcaptionbox{}{\includegraphics[width=0.45\linewidth]{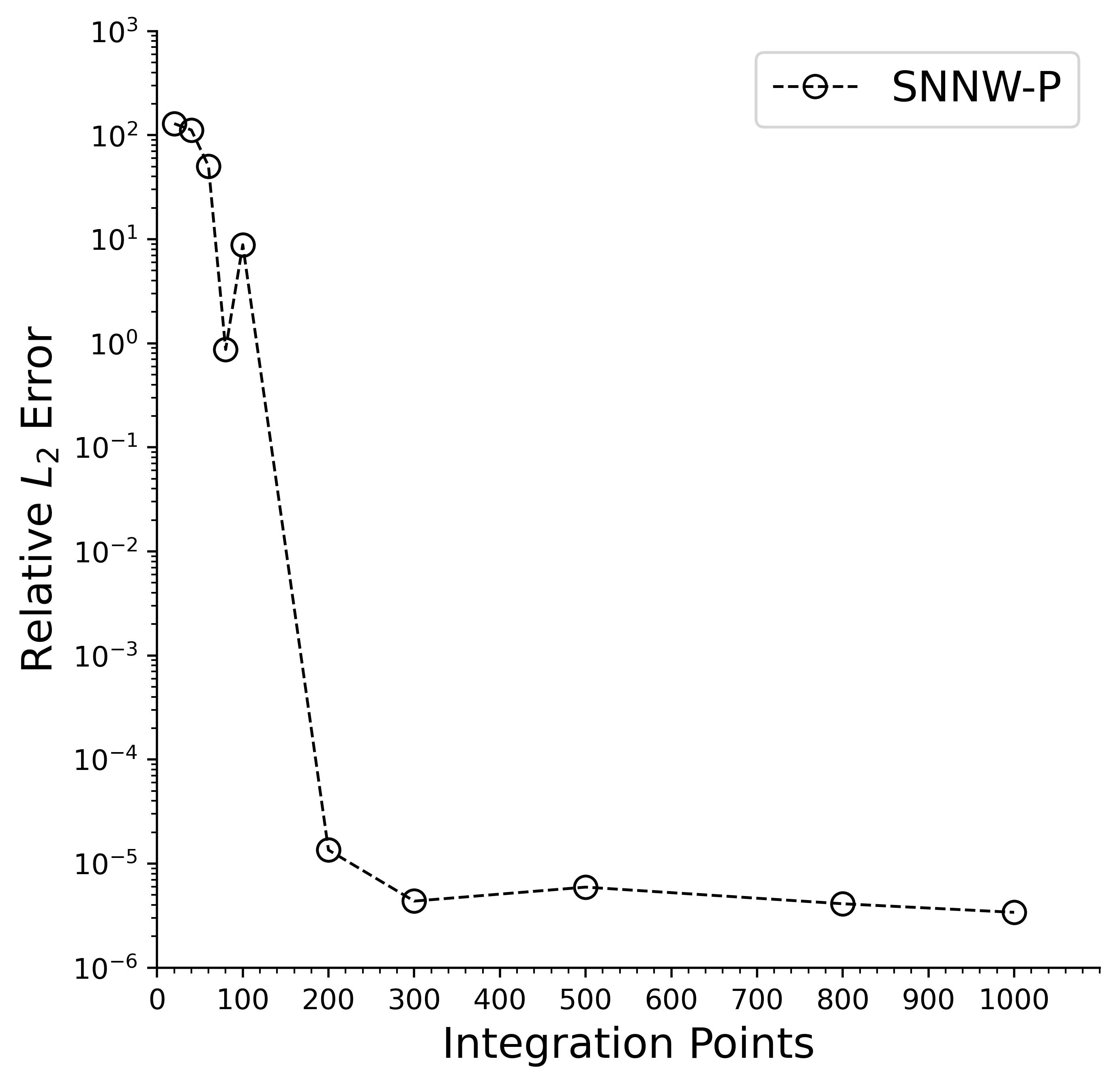}}
\caption{Error variation with subspace dimension at a fixed number of 1000 sampling points and error variation with the number of sampling points at a fixed subspace dimension of 300 for SNNW-P.}
\label{hz_SNNW_P_subspace_fig}
\end{figure}

\begin{figure}[htbp]
\centering
\subcaptionbox{}{\includegraphics[width=0.45\linewidth]{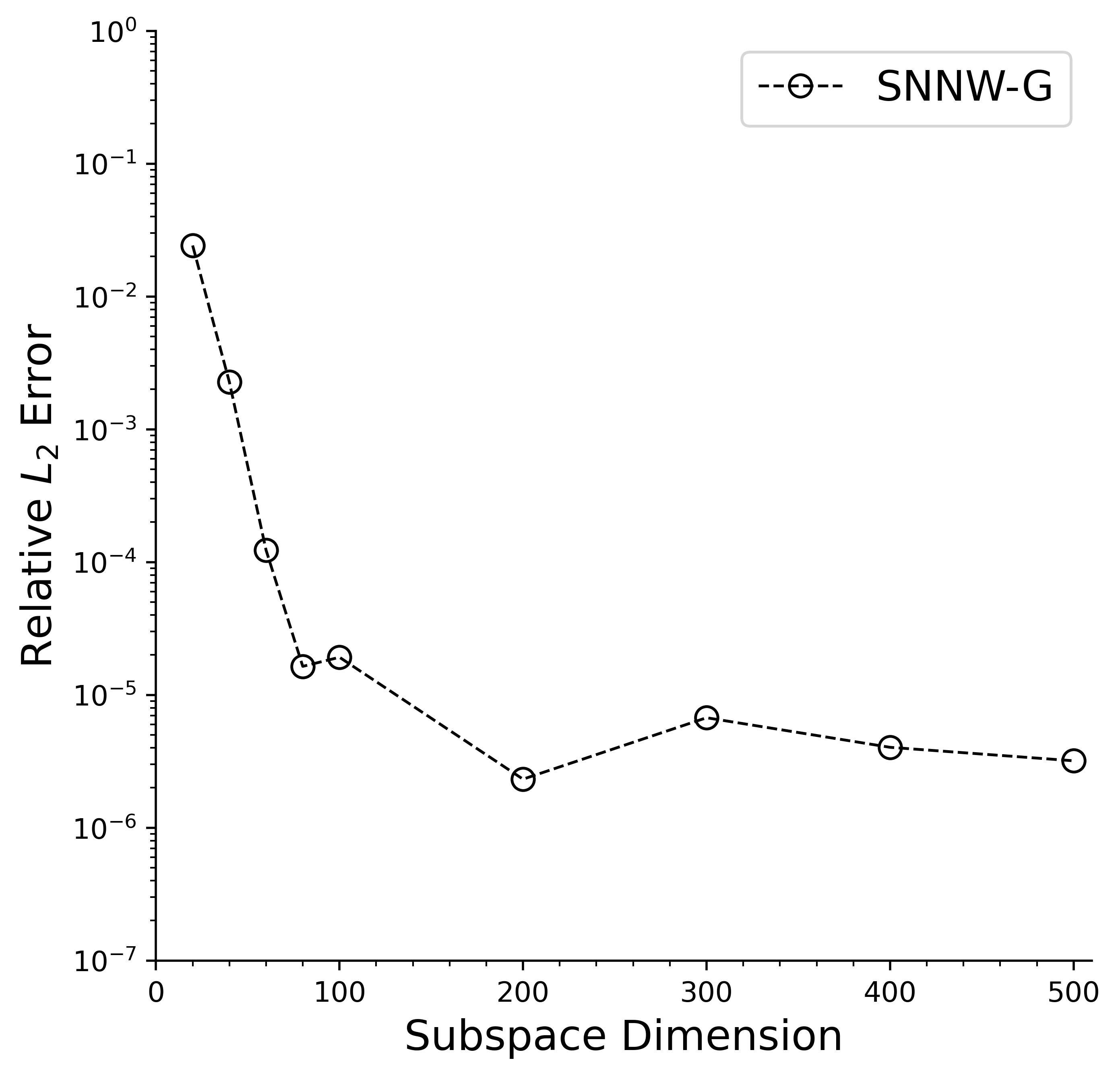}}
\hspace{0.3cm}  
\subcaptionbox{}{\includegraphics[width=0.45\linewidth]{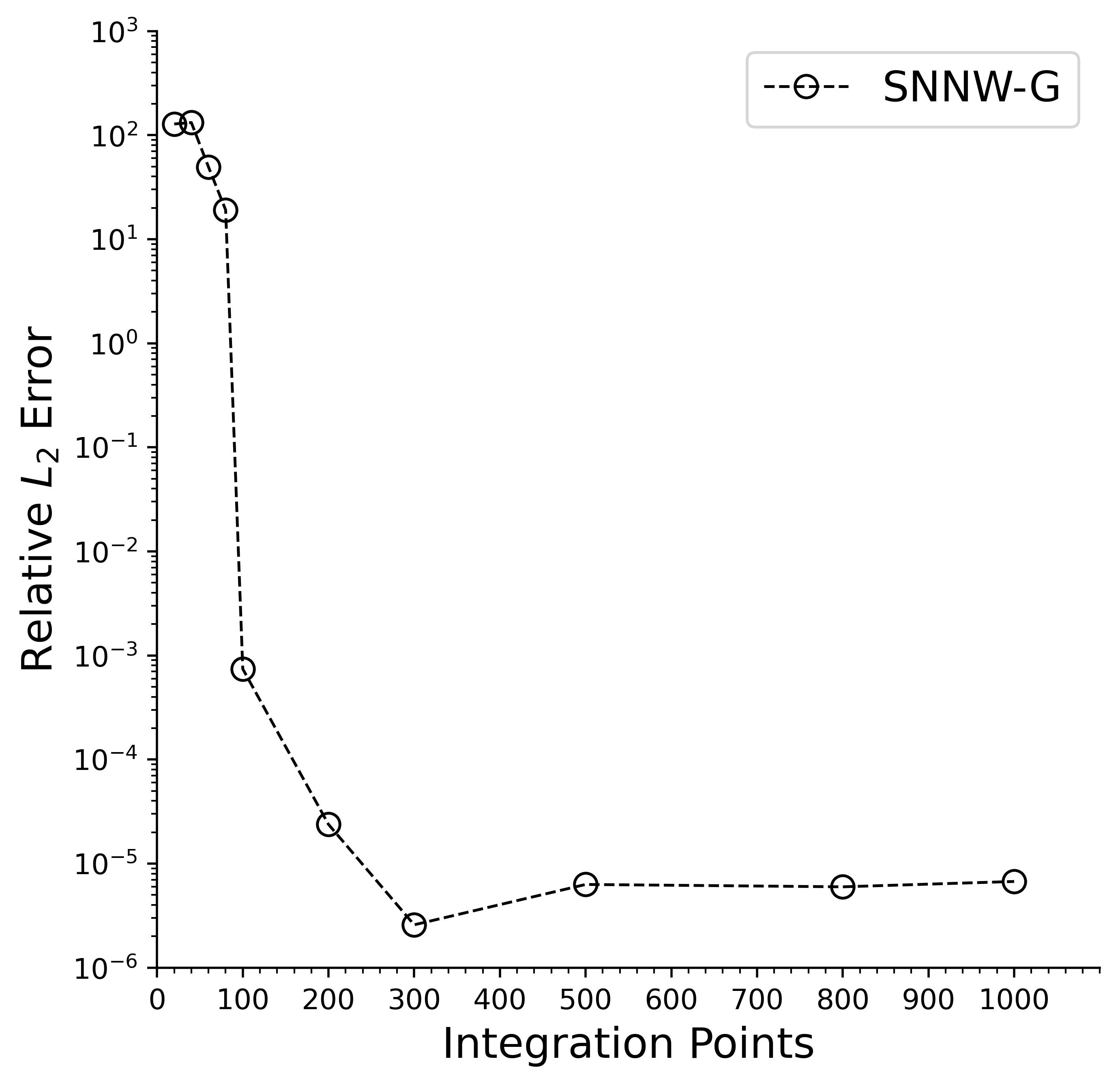}}
\caption{Error variation with subspace dimension at a fixed number of 1000 sampling points and error variation with the number of sampling points at a fixed subspace dimension of 300 for SNNW-G.}
\label{hz_SNNW_G_subspace_fig}
\end{figure}

\begin{figure}[htbp]
\centering
\subcaptionbox{}{\includegraphics[width=0.45\linewidth]{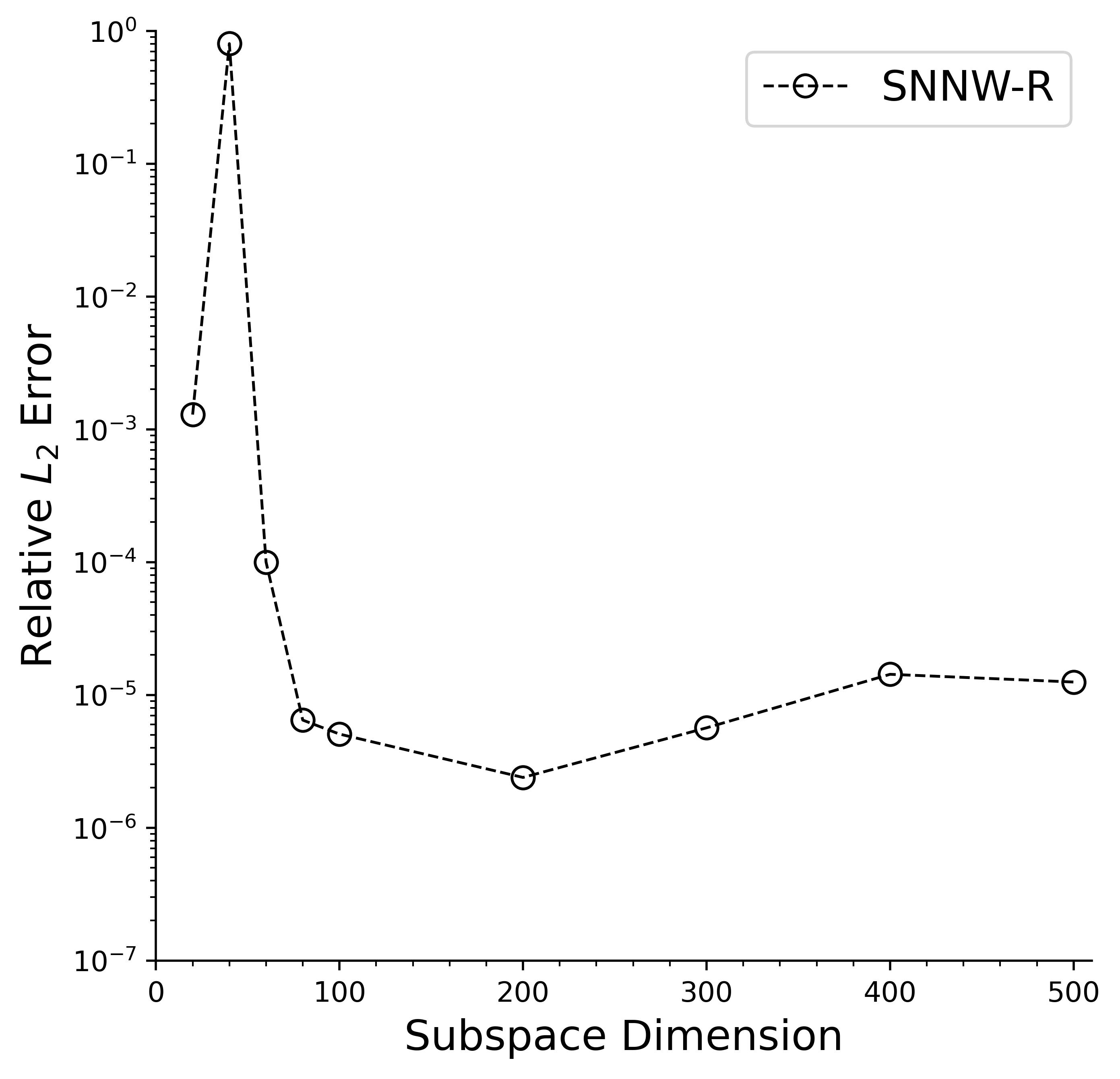}}
\hspace{0.3cm}  
\subcaptionbox{}{\includegraphics[width=0.45\linewidth]{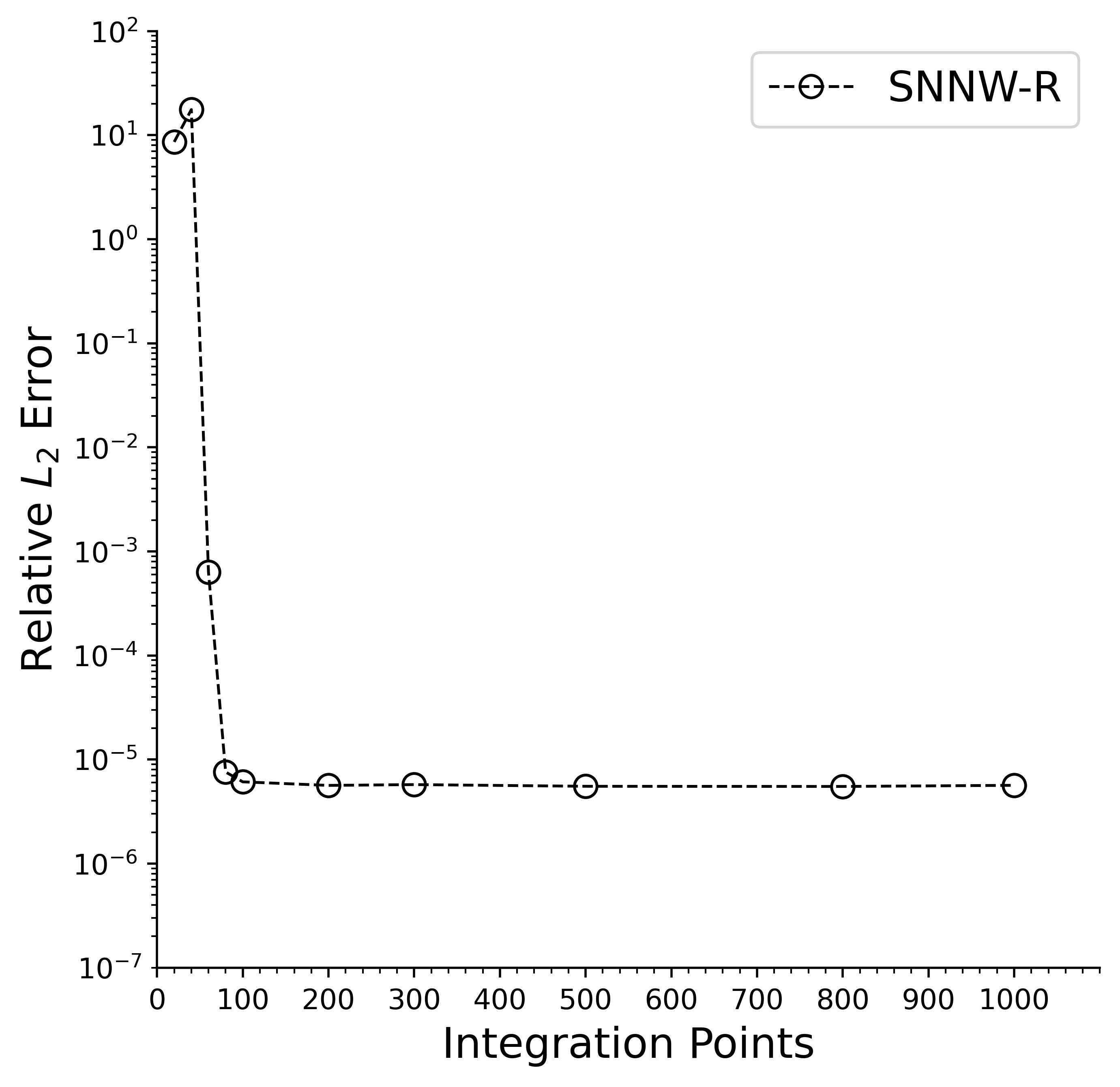}}
\caption{Error variation with subspace dimension at a fixed number of 1000 sampling points and error variation with the number of sampling points at a fixed subspace dimension of 300 for SNNW-R.}
\label{hz_SNNW_R_subspace_fig}
\end{figure}

Now we test the variation pattern of errors. First, the number of sampling points is given, we examine the error variation with respect to subspace dimension.
Then,  the subspace dimension is fixed, we examine the error variation with respect to
the number of sampling points. Figure \ref{hz_SNNW_P_subspace_fig} illustrates the error variation with subspace dimension for 1000 sampling points, and the error variation with
the number of sampling points for a fixed subspace dimension of 300 for
SNNW-P. Similarly, Figure  \ref{hz_SNNW_G_subspace_fig} and  \ref{hz_SNNW_R_subspace_fig} show the variation pattern of errors for SNNW-G and SNNW-R, respectively.

Table \ref{hz_SNN_D_netsize_subspace} presents the results of SNNW-P for different hidden layer depth and subspace dimension, maintaining a uniform width of 100 for each hidden layer.
When the number of the  subspace dimension is 60, the accuracy is not satisfied.
When the number of the  subspace dimension is 300, the error of SNNW-P is between $10^{-5}$ and $10^{-7}$, and the depth of hidden layer has no   significant impact on accuracy.
For SNNW-G and SNNW-R, there has similar results.

\begin{table}[!htbp]
\caption{The errors and epochs for SNNW-P on different hidden layer depth.}\label{hz_SNN_D_netsize_subspace}
\begin{center}
\small
\begin{tabular}{ccccc}
\hline
Hidden Layer & M  & 60 & 100 & 300 \\
\hline
0 & $\|e\|_{ {L}^2}$ & 1.11e-05 & 1.27e-05 & 2.39e-06 \\
& epochs & 3340& 2600 &  1902 \\
1 & $\|e\|_{ {L}^2}$ &  2.01e-03 &  3.39e-04 & 3.33e-07 \\
& epochs &   2540& 860& 254 \\
2 & $\|e\|_{ {L}^2}$ & 3.49e-01 & 8.29e-05 &  2.42e-05\\
& epochs &  5000&   1415& 406 \\
3 & $\|e\|_{ {L}^2}$ & 2.41e-04 & 1.42e-05 & 1.94e-06 \\
& epochs & 393 &  233 & 213\\
4 & $\|e\|_{ {L}^2}$ & 6.67e-05 &  2.714e-05& 3.39e-06\\
& epochs &  295& 676 & 282\\
5 & $\|e\|_{ {L}^2}$ & 1.01e-03& 1.32e-04& 3.90e-06\\
& epochs & 1465 & 456& 275 \\
6 & $\|e\|_{ {L}^2}$ &5.57e-04 & 1.14e-04&  1.45e-06\\
& epochs & 4159&  717&190\\
\hline
\end{tabular}
\end{center}
\end{table}

\begin{table}[!htbp]
\caption{The errors and epochs for SNNW-G on different hidden layer depth.}\label{hz_SNNW_G_netsize_subspace}
\begin{center}
\small
\begin{tabular}{ccccc}
\hline
Hidden Layer & M  & 60 & 100 & 300 \\
\hline
0 & $\|e\|_{ {L}^2}$ & 1.76e-04 & 9.40e-06 & 1.30e-06 \\
& epochs & 4026& 2488 &1894   \\
1 & $\|e\|_{ {L}^2}$ &   2.61e-03& 2.99e-04  & 2.10e-07 \\
& epochs &   2455 & 912&  215 \\
2 & $\|e\|_{ {L}^2}$ & 2.94e-02& 2.59e-04 & 4.81e-05  \\
& epochs & 4550 &  1356 &  360 \\
3 & $\|e\|_{ {L}^2}$ & 1.44e-04 & 3.64e-05 & 1.80e-06\\
& epochs &  380 & 349 &196\\
4 & $\|e\|_{ {L}^2}$ &  1.23e-04&  1.92e-05& 6.73e-06\\
& epochs & 314&  362&  273\\
5 & $\|e\|_{ {L}^2}$ & 6.90e-04& 5.98e-06& 6.32e-06\\
& epochs &  757& 335& 208 \\
6 & $\|e\|_{ {L}^2}$ &2.25e-04 & 2.00e-04 & 1.50e-05 \\
& epochs & 1717&  553 & 307\\
\hline
\end{tabular}
\end{center}
\end{table}

\begin{table}[!htbp]
\caption{The errors   for SNNW-R on different hidden layer depth.}\label{hz_SNNW_R_netsize_subspace}
\begin{center}
\small
\begin{tabular}{ccccc}
\hline
Hidden Layer & M  & 60 & 100 & 300 \\
\hline
0 & $\|e\|_{ {L}^2}$ & 8.40e-05  & 2.88e-05 & 3.83e-06 \\
1 & $\|e\|_{ {L}^2}$ & 9.35e-04  &  8.91e-06 &  3.54e-07\\
2 & $\|e\|_{ {L}^2}$ & 3.71e-04& 4.14e-05 &  4.30e-03 \\
3 & $\|e\|_{ {L}^2}$ & 8.97e-05 &  9.09e-07& 1.15e-05 \\
4 & $\|e\|_{ {L}^2}$ & 1.00e-04& 5.07e-06 &5.65e-06 \\
5 & $\|e\|_{ {L}^2}$ &2.36e-04 & 2.94e-05& 4.88e-06\\
6 & $\|e\|_{ {L}^2}$ & 2.91e-04&  1.09e-05&  5.16e-06\\
\hline
\end{tabular}
\end{center}
\end{table}

\newpage
\subsection{Poisson equation} \label{Poisson equation}
Consider the following Poisson equation with the Dirichlet boundary condition,
\begin{eqnarray*}
-\Delta u(x,y)= f,\,(x,y)\in\Omega,
\end{eqnarray*}
where $ \Omega=(0,1)^2$. The Dirichlet boundary condition is
\begin{eqnarray*}
u(x,y)=0,\,(x,y)\in\partial\Omega.
\end{eqnarray*}
Take the exact solution as follows:
\begin{eqnarray*}
u(x,y)=\sin{(\pi x)}\sin{(\pi y)}.
\end{eqnarray*}

\begin{figure}[htbp]
\centering
\subcaptionbox{SNNW-P}{\includegraphics[width=0.3\linewidth]{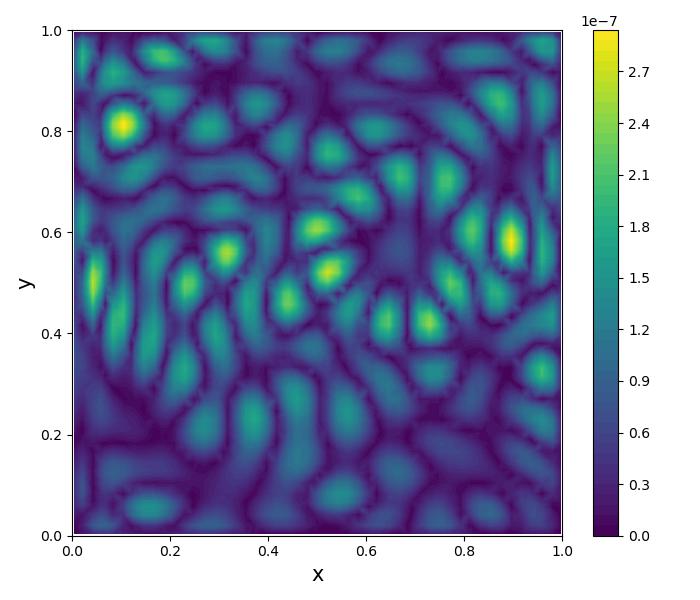}}
\hspace{0.3cm}  
\subcaptionbox{SNNW-G}{\includegraphics[width=0.3\linewidth]{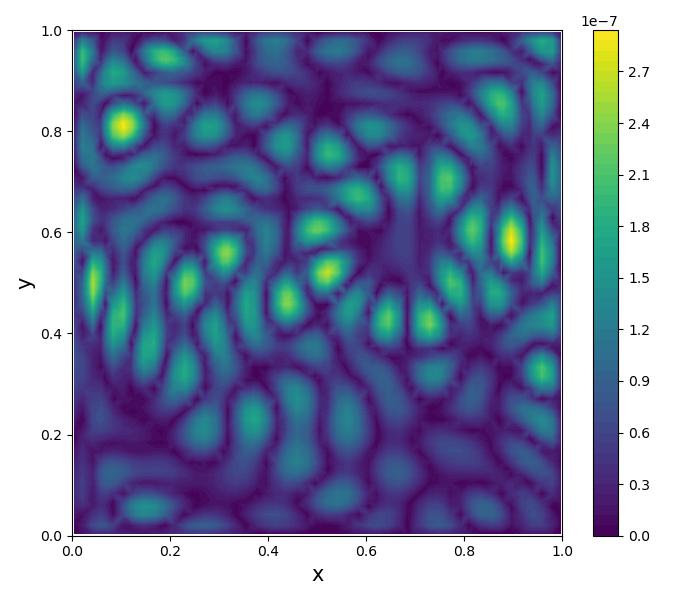}}
\hspace{0.3cm} 
\subcaptionbox{SNNW-R}{\includegraphics[width=0.3\linewidth]{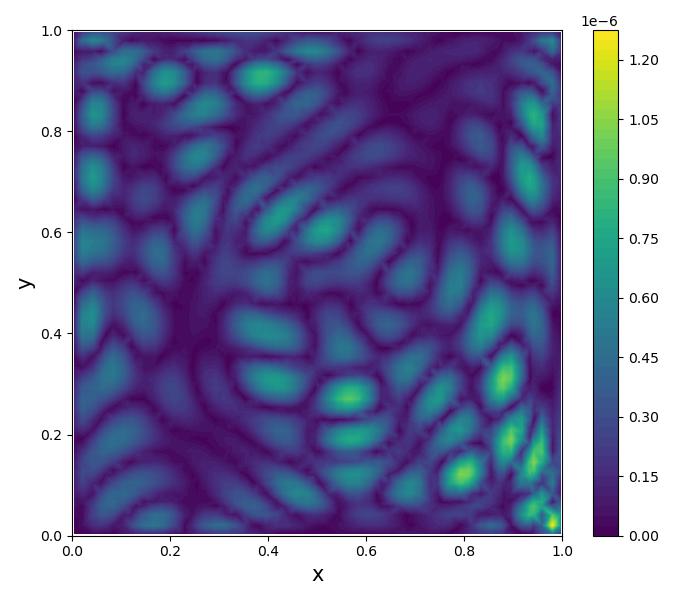}}
\caption{ Point-wise errors  of SNNW-P, SNNW-G and SNNW-R for Poisson equation.}
\label{ps_point_error}
\end{figure}

For SNNW, a two-dimensional composite Gaussian quadrature formula is applied to segment each dimension into 16 sub-intervals and allocate 4 points to
each.
The point-wise errors of  SNNW-P, SNNW-G and SNNW-R  are illustrated in Figure \ref{ps_point_error}.
Table \ref{possion different method} presents the relative $L^2$ errors for various methods including SNNW-P, SNNW-G, SNNW-R, PINN, DGM and DRM. This example serves to
demonstrate the adaptivity of our algorithms. Notably, after 50000 training
epochs, the relative $L^2$  errors of PINN, DGM and DRM are 2.71e-04, 1.81e-03 and 2.48e-04,
respectively. However, SNNW-P, SNNW-G and SNNW-R achieve relative $L^2$  errors of 1.67e-07, 5.79e-07
 and 1.65e-07, respectively, with significantly fewer epochs, specifically  359, 359 and 2000.

\begin{table}[!htbp]
\caption{The errors and epochs of different methods for the Poisson equation}\label{possion different method}
\begin{center}
\small
\begin{tabular}{cccc}\hline
Method & $\|e\|_{L^2}$ & epochs\\
\hline
PINN & 2.71e-04 & 50000 \\
DGM &  1.81e-03 & 50000 \\
DRM & 2.48e-04 &  50000\\
SNNW-P & 1.67e-07 & 359 \\
SNNW-G & 1.65e-07 & 359 \\
SNNW-R & 5.79e-07 & 2000\\
\hline
\end{tabular}
\end{center}
\end{table}

Tables \ref{possionSNNWP} and \ref{possionSNNWG} present the performance of SNNW-P and SNNW-G across
various numbers of sampling points and subspace dimension, maintaining a
networks architecture of four hidden layers with 100 neurons each. SNNW-P and SNNW-G adopt the composite Gaussian
quadrature, allocating 4 points per sub-interval in each direction and increasing sub-interval numbers for enhanced precision.
Similar to the 1D Helmholtz equation, it becomes apparent that the error decreases with increasing both sampling points and subspace dimension.
With a subspace dimension of 40, the error primarily ranges between $10^{-1}$ and $10^{-4}$.
Expanding the subspace dimension to 100, SNNW-P and SNNW-G can achieve the accuracy of $10^{-5}$.
With enough sampling points and further expansion of the subspace dimension, SNNW-P and SNNW-G can obtain and maintain the  accuracy of $10^{-7}$.
Moreover, the training cost is low, and there only needs  300 to 2000 epochs.

\begin{table}[!htbp]
\caption{The errors and epochs of SNNW-P across various numbers of sampling points and subspace dimension $M$ for the Poisson equation}\label{possionSNNWP}
\begin{center}
\resizebox{\columnwidth}{!}{%
\small
\begin{tabular}{ccccccccc}
\hline
Points & $M$ & 40 & 60 & 80 & 100 & 200 & 300\\
\hline
$12\times12$ & $\|e\|_{ {L}^2}$ & 2.26e-01 & 1.24e-01 & 1.22e-01 & 7.17e-02 & 6.66e-01 & 7.75e-01 \\
 & epochs  & 1482 & 936 & 669 & 572 & 398 & 320 \\
$16\times16$ & $\|e\|_{ {L}^2}$& 7.18e-03 & 7.66e-03 & 1.38e-02 & 2.34e-03 & 7.30e-02 & 1.4e-01  \\
& epochs & 1536 & 972 & 701 & 612 & 427 & 340 \\
$20\times20$ & $\|e\|_{L^2}$ & 1.71e-03 & 2.18e-03 & 6.25e-04 & 2.72e-04 & 1.05e-02 & 1.13e-02\\
& epochs & 1452 & 1009 & 733 & 635 & 439 & 340 \\
$24\times24$ & $\|e\|_{L^2}$ & 2.48e-04 & 2.42e-04 & 9.55e-05 & 8.20e-05 & 7.64e-04 & 8.00e-04 \\
&epochs& 1440 & 1029 & 743 & 649 & 445 & 354 \\
$32\times32$&$ \|e\|_{L^2} $ & 1.97e-04 & 1.45e-04 & 6.80e-05 & 2.11e-05 & 4.48e-05 & 4.29e-05\\
&epochs& 1388 & 1069 & 753 & 660 & 454 & 357 \\
$48\times48$ & $\|e\|_{L^2}$ & 2.36e-04 & 1.27e-04 & 4.23e-05 & 1.48e-05 & 5.74e-07 & 3.89e-07 \\
& epochs& 1502 & 1076 & 768 & 674 & 459 & 361\\
$64\times64$ & $\|e\|_{L^2}$ & 3.07e-04 & 1.45e-04 & 7.56e-05 & 1.41e-05 & 2.06e-07 & 1.67e-07 \\
&epochs & 1392 & 1066 & 767 & 677 & 461 & 359 \\
\hline
\end{tabular}
} 
\end{center}
\end{table}

\begin{table}[!htbp]
\caption{The errors and epochs of SNNW-G across various numbers of sampling points and subspace dimension $M$ for the Poisson equation}\label{possionSNNWG}
\begin{center}
\resizebox{\columnwidth}{!}{%
\small
\begin{tabular}{ccccccccc}
\hline
Points & $M$ & 40 & 60 & 80 & 100 & 200 & 300\\
\hline
$12\times12$ & $\|e\|_{ {L}^2}$ & 4.16e-02  & 1.26e-01  &  2.91e-01 & 7.11e-02 & 6.86e-01 & 7.71e-01  \\
 & epochs  & 1846 &947 & 668 & 571 & 398 & 320\\
$16\times16$ & $\|e\|_{ {L}^2}$& 5.39e-03  & 9.60e-04 & 1.00e-02 &  2.29e-03 & 7.22e-01 & 1.42e-01  \\
& epochs & 1527 & 988 & 716 & 612 & 427 & 340\\
$20\times20$ & $\|e\|_{L^2}$ & 3.31e-04 & 3.29e-04 & 6.32e-04 & 2.70e-04 & 1.02e-02 & 1.21e-02\\
& epochs & 1380  & 1010 & 733 & 635  & 439  & 348 \\
$24\times24$ & $\|e\|_{L^2}$ & 2.02e-04 & 3.17e-04 & 3.95e-04 & 8.31e-05 & 6.03e-04 & 7.99e-04 \\
&epochs& 1380  & 1035 & 741 &  650 &  445 & 354 \\
$32\times32$&$ \|e\|_{L^2} $ & 2.52e-04 & 1.41e-04 & 6.64e-05 &  2.14e-05 & 4.60e-05 & 4.31e-05\\
&epochs& 1487 &  1070 &  752 & 659 & 454 &  357\\
$48\times48$ & $\|e\|_{L^2}$ & 2.61e-04 & 1.12e-04 & 3.81e-05 & 1.53e-05 & 5.87e-07  & 3.83e-07 \\
& epochs& 1451 & 1065 & 770 & 675 &  459 & 361 \\
$64\times64$ & $\|e\|_{L^2}$ & 2.42e-04  & 1.45e-04 & 6.26e-05 & 1.44e-05 & 2.07e-07  &  1.65e-07  \\
&epochs & 1426  & 1070 & 768 & 676  & 461 & 359 \\
\hline
\end{tabular}
} 
\end{center}
\end{table}

The network setup and composite Gaussian quadrature rule of SNNW-R are identical with the ones of SNNW-P and SNNW-G.
The numerical results of SNNW-R with fixed 2000 epochs are given in Table \ref{possionSNNWR} for various numbers
of sampling points and subspace dimension. With a subspace dimension of 100 and 150, the error is about $10^{-3}$. Increasing the subspace dimension to 200, the  error decreases promptly to $10^{-6}$.
As the number of sampling points and  the subspace dimension further increase, SNNW-R achieves the accuracy of  $10^{-7}$.

\begin{table}[!htbp]
\caption{The errors of SNNW-R with 2000 epochs across various numbers of sampling points and subspace dimension $M$ for the Poisson equation}\label{possionSNNWR}
\begin{center}
\resizebox{\columnwidth}{!}{%
\small
\begin{tabular}{ccccccccc}
\hline
Points & $M$ & 100 & 150 & 200 & 250 & 300 & 500\\
\hline
$20\times20$ & $\|e\|_{L^2}$ & 1.42e-03 & 5.05e-03 & 2.24e-02 & 5.88e-01  & 5.53e-01  & 3.00e-01 \\
$24\times24$ & $\|e\|_{L^2}$ & 1.42e-03 & 2.21e-03 & 1.60e-03 & 1.35e-04  & 1.73e-03  & 6.56e-03 \\
$32\times32$ & $\|e\|_{L^2}$ & 1.97e-03 & 1.76e-03 & 1.15e-04 & 2.26e-04  & 1.06e-04  & 2.34e-04 \\
$48\times48$ & $\|e\|_{L^2}$ & 1.91e-03 & 1.66e-03 & 3.23e-06 & 4.26e-06  & 3.23e-06  & 1.21e-05 \\
$64\times64$ & $\|e\|_{L^2}$ & 1.65e-03 & 1.67e-03 & 2.12e-06 & 8.00e-07  & 5.79e-07  & 9.32e-07 \\
$72\times72$ & $\|e\|_{L^2}$ & 2.05e-03 & 1.72e-03 & 2.10e-06 & 7.15e-07  & 5.18e-07  & 3.31e-07  \\
\hline
\end{tabular}
} 
\end{center}
\end{table}

Figure 7 illustrates the error variation with subspace dimension for $64\times64$
Gaussian integration points, and the error variation with respect to the number of  integration points
for a fixed subspace dimension of 300 for SNNW-P.
Similarly, Figure 8 and 9 show the error
variation with subspace dimension at a fixed number of $64\times64$
Gaussian integration points, and the error variation with the number of Gaussian integration points at a fixed subspace dimension of 300 for SNNW-G and SNNW-R, respectively.  The number of training epochs
for these experiments is ranging from 359 to 2000, and is significantly lower than that of PINN, DGM
and DRM.
\begin{figure}[htbp]
\centering
\subcaptionbox{}{\includegraphics[width=0.45\linewidth]{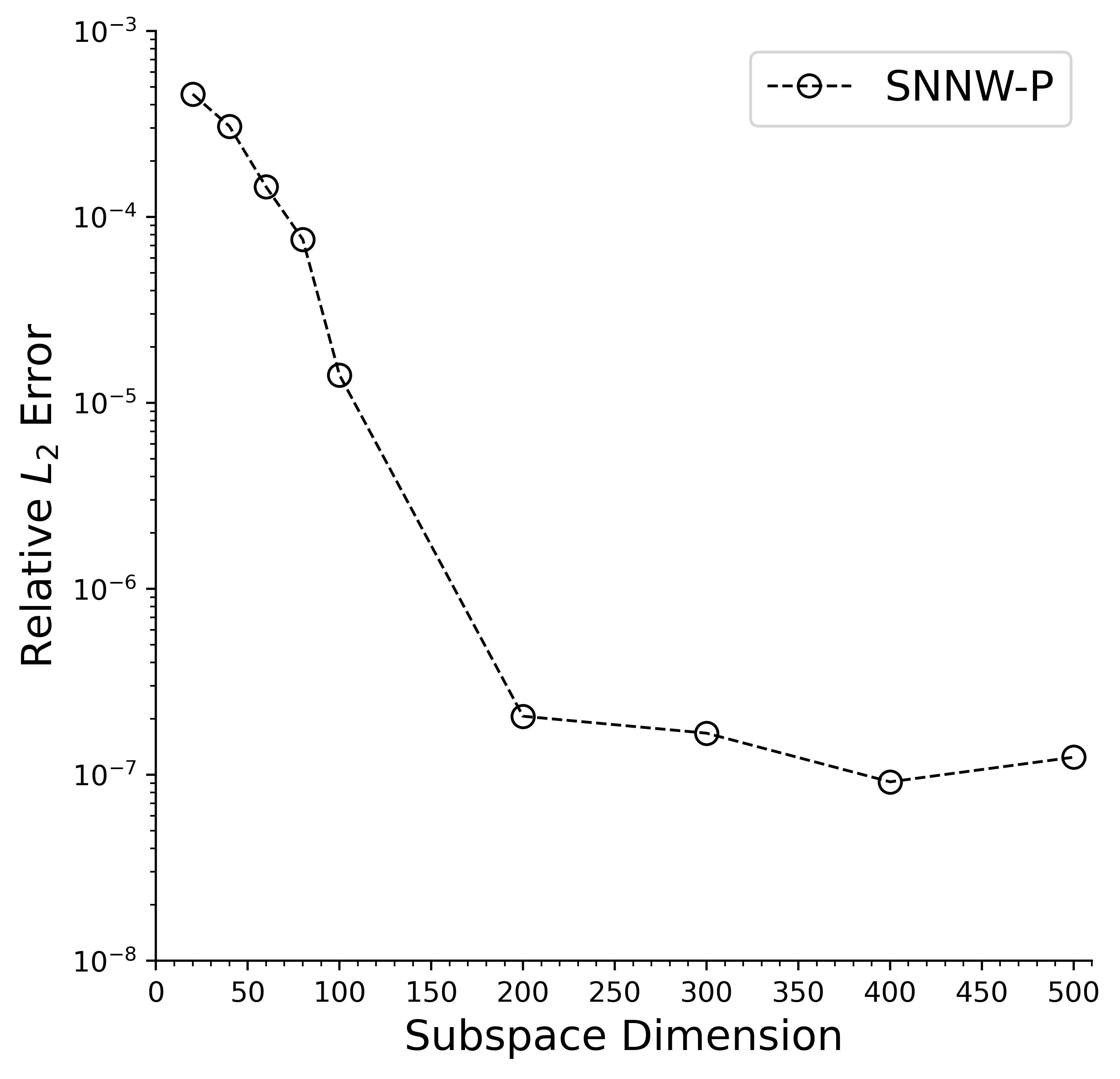}}
\hspace{0.3cm}  
\subcaptionbox{}{\includegraphics[width=0.45\linewidth]{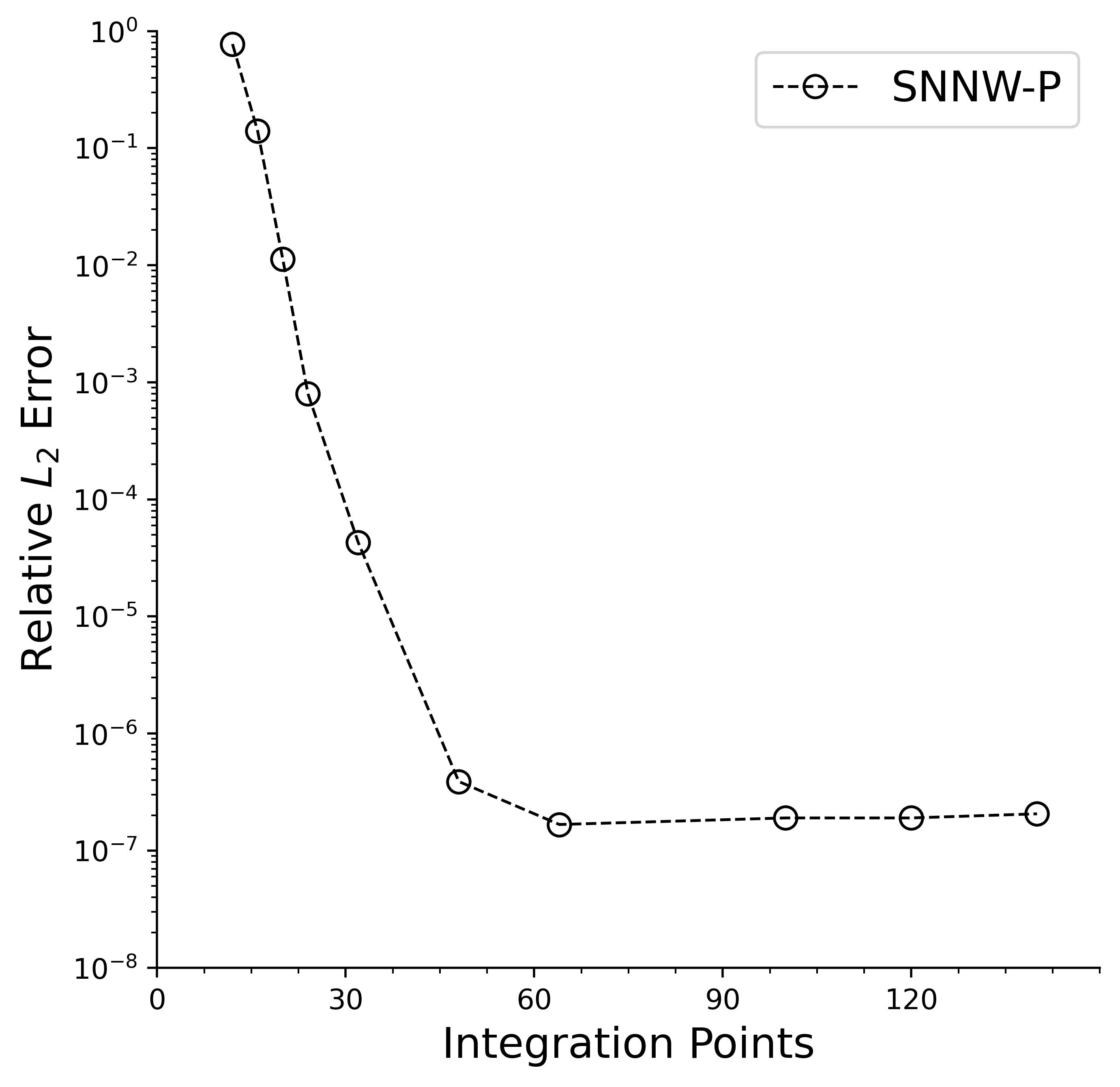}}
\caption{Error variation with subspace dimension at a fixed number of \(64 \times 64\) integration points and error variation with the number of integration points at a fixed subspace dimension of 300 for SNNW-P.}
\label{2dps_SNN_I_subspace_fig}
\end{figure}

\begin{figure}[htbp]
\centering
\subcaptionbox{}{\includegraphics[width=0.45\linewidth]{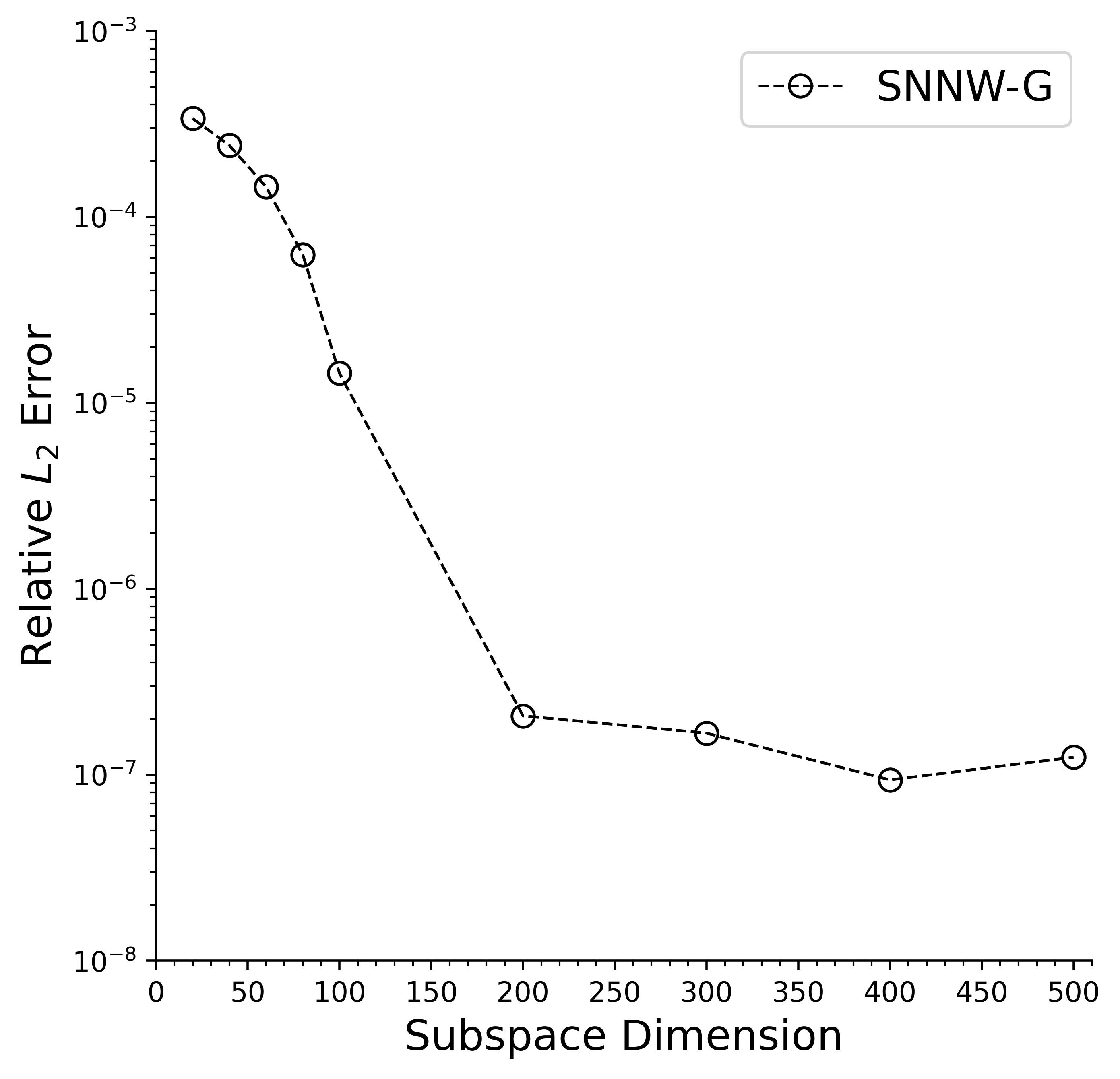}}
\hspace{0.3cm}  
\subcaptionbox{}{\includegraphics[width=0.45\linewidth]{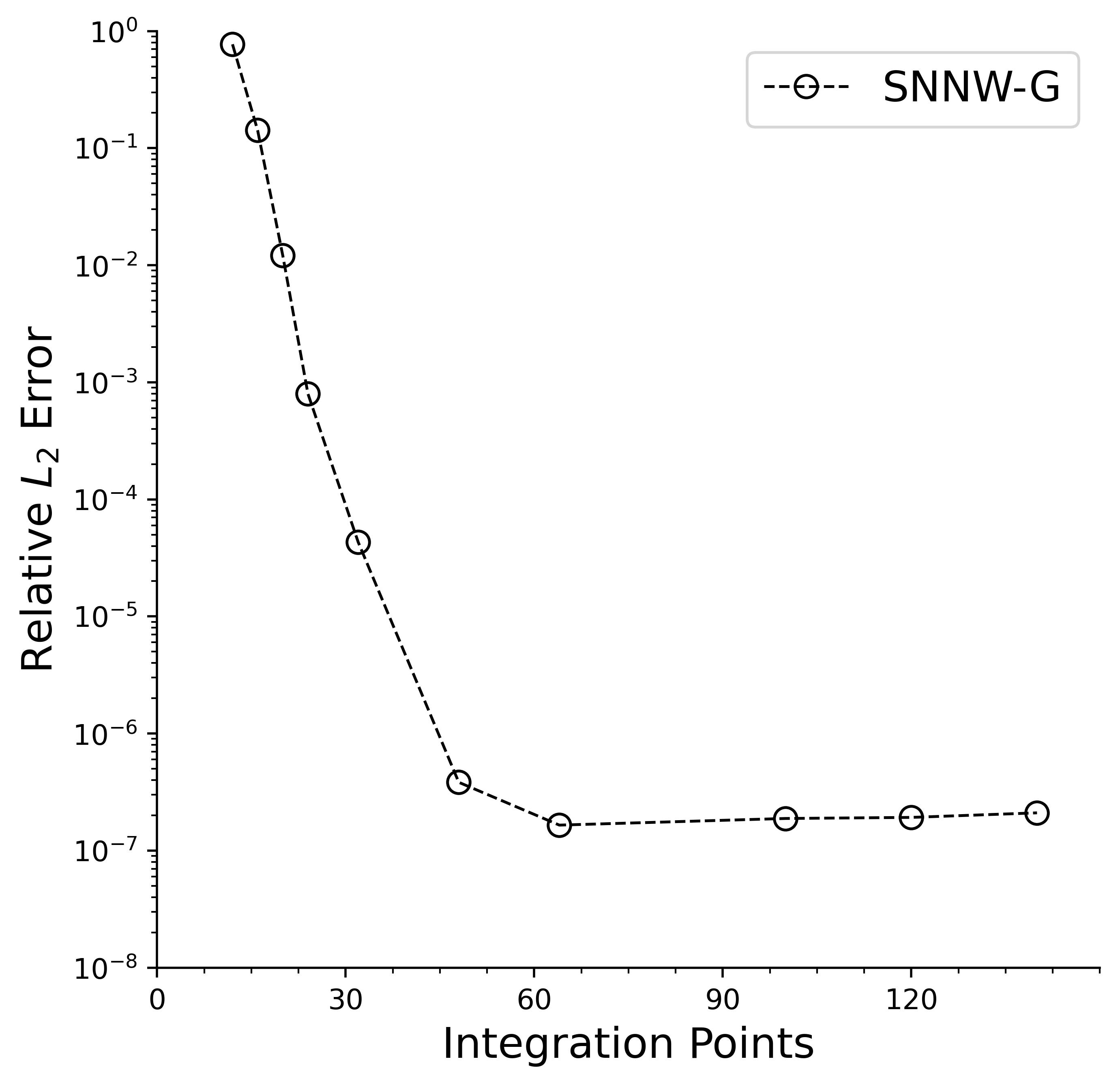}}
\caption{Error variation with subspace dimension at a fixed number of \(64 \times 64\) integration points and error variation with the number of integration points at a fixed subspace dimension of 300 for SNNW-G.}
\label{2dps_SNN_I_subspace_fig}
\end{figure}

\begin{figure}[htbp]
\centering
\subcaptionbox{}{\includegraphics[width=0.45\linewidth]{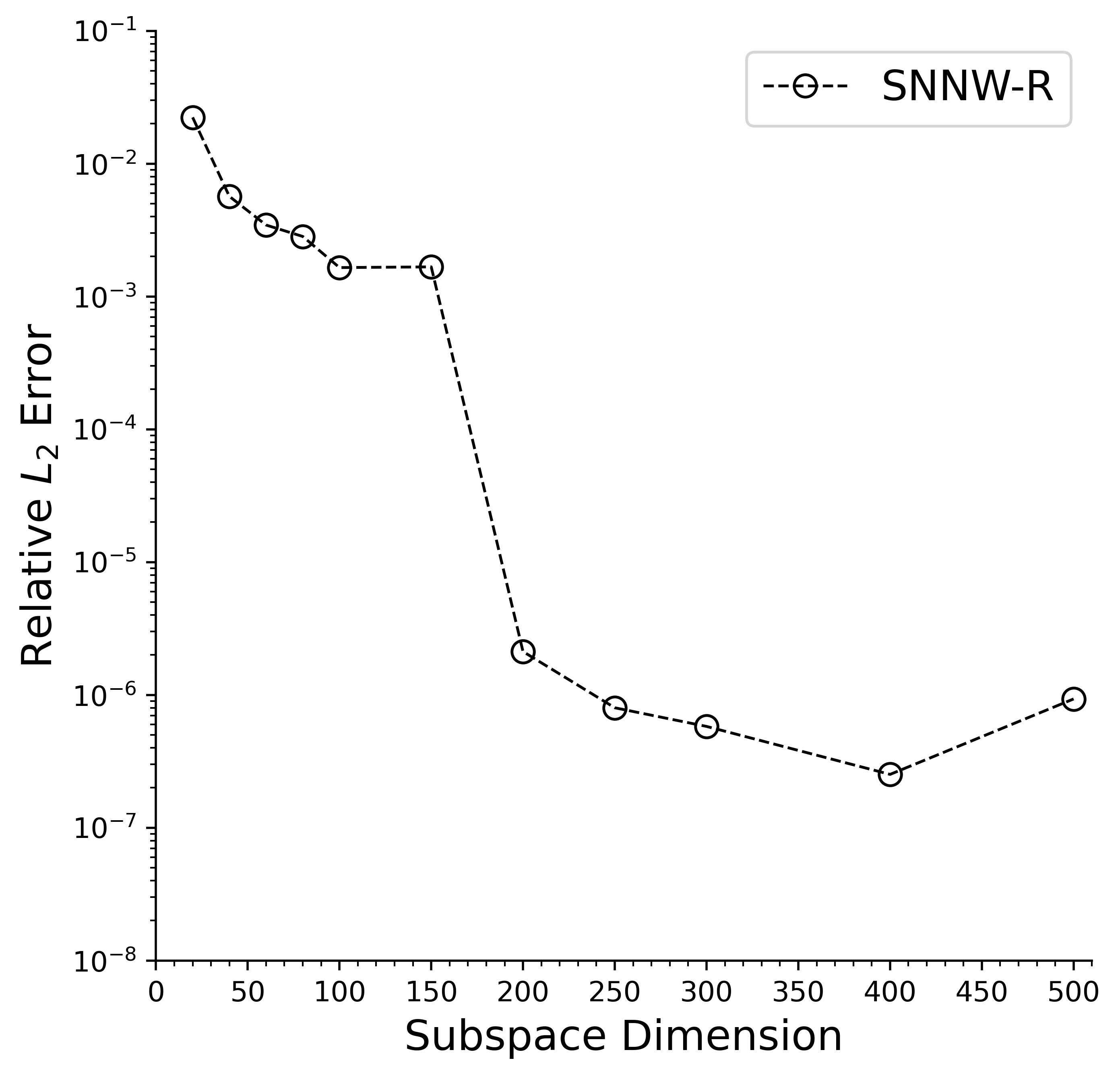}}
\hspace{0.3cm}  
\subcaptionbox{}{\includegraphics[width=0.45\linewidth]{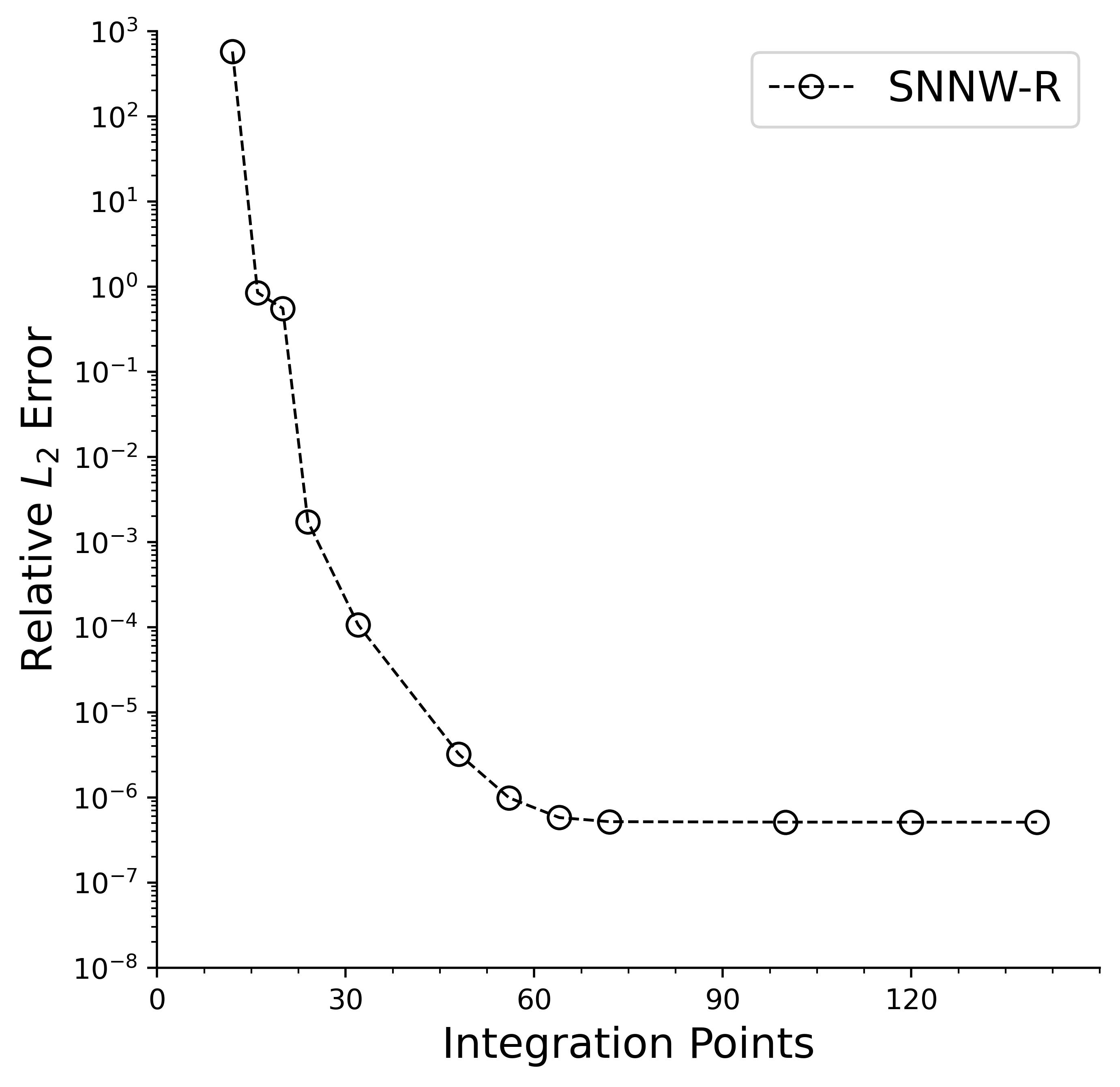}}
\caption{Error variation with subspace dimension at a fixed number of \(64 \times 64\) integration points and error variation with the number of integration points at a fixed subspace dimension of 300 for SNNW-R.}
\label{2dps_SNN_I_subspace_fig}
\end{figure}

\newpage
\subsection{Anisotropic diffusion equation}\label{Anisotropic diffusion equation}
Consider the following  diffusion equation with anisotropic diffusion coefficient on the domain \(\Omega=(0,1)^2\),
\begin{eqnarray*}
    \label{Anisotropic_eq}
    -\nabla \cdot(\kappa(x, y) \nabla u)=f(x, y), & (x, y) \in \Omega, \\
    u(x, y)=g(x, y), & (x, y) \in \partial \Omega,
\end{eqnarray*}
where
\begin{eqnarray*}
\kappa(x, y)=\left(\begin{array}{cc}
k_1 & 0 \\
0 & k_2
\end{array}\right).
\end{eqnarray*}
We choose  $f (x,y)$ and $g (x,y)$ such that the exact solution is
$
u(x,y) = \sin\left(\pi x\right) \sin\left(\pi y \right).
$

A two-dimensional composite Gaussian quadrature formula is utilized to divide each dimension into 16 sub-intervals, with 4 points per sub-interval.
This strategy leads to 4096 sampling points.
Table \ref{ADclassicmethod} and  \ref{ADSNNW} illustrate the relative $L^2$ errors of different methods for solving the anisotropic diffusion equation with various anisotropic ratios.
PINN, DGM and DRM can  persevere certain accuracy when the anisotropy is  weak, and lose the accuracy when the anisotropy is very strong.
However,  even when the anisotropic ratio reaches $1:10^6$, SNNW-P attains the accuracy of 2.02e-06, SNNW-G attains the accuracy of 2.07e-06 and SNNW-R attains the accuracy of 1.56e-05. This indicates that SNNW-P, SNNW-G and SNNW-R still can obtain high accuracy for the diffusion equation with strong anisotropic ratios.
Compared to  50000 epochs of PINN, DGM and DRM, SNNW-P and SNNW-G only need several hundred of epochs.
These tests further demonstrate that our method can achieves  better accuracy with low training cost.

\begin{table}[!htbp]
\caption{The errors and epochs of PINN, DGM and DRM  for the anisotropic diffusion equation with different anisotropic ratios.}\label{ADclassicmethod}
\begin{center}
\small
\begin{tabular}{ccccc}
\hline
Method & $k_1:k_2$    & $\|e\|_{L^2}$            & epochs \\
\hline
       & $1:10^0$     &   9.54e-03    & 50000      \\
       & $1:10^1$     &   3.14e-02    & 50000      \\
       & $1:10^2$     &   1.68e-01    & 50000      \\
PINN   & $1:10^3$     &   1.68e00     & 50000      \\
       & $1:10^4$     &   1.03e01     & 50000      \\
       & $1:10^5$     &   9.62e01     & 50000      \\
       & $1:10^6$     &   1.04e01     & 50000      \\
\hline
       & $1:10^0$     &  1.81e-03     & 50000      \\
       & $1:10^1$     &  2.67e-02     & 50000      \\
       & $1:10^2$     &  7.55e-01     & 50000      \\
 DGM   & $1:10^3$     &  1.91e00      & 50000      \\
       & $1:10^4$     &  1.00e01      & 50000      \\
       & $1:10^5$     &  1.02e01      & 50000      \\
       & $1:10^6$     &  1.02e01      & 50000      \\
\hline
       & $1:10^0$     &   2.48e-04    & 50000      \\
       & $1:10^1$     &   1.44e-04    & 50000      \\
       & $1:10^2$     &   1.37e-04    & 50000      \\
 DRM   & $1:10^3$     &   2.16e-01    & 50000      \\
       & $1:10^4$     &   4.77e-02    & 50000      \\
       & $1:10^5$     &   4.97e-04    & 50000      \\
       & $1:10^6$     &   2.06e-02    & 50000      \\
\hline
\end{tabular}
\end{center}
\end{table}

\begin{table}[!htbp]
\caption{The errors and epochs of SNNW-P,SNNW-G and SNNW-R for the anisotropic diffusion equation with different anisotropy ratios.}\label{ADSNNW}
\begin{center}
\small
\begin{tabular}{ccccc}
\hline
Method & $k_1:k_2$    & $\|e\|_{L^2}$            & epochs \\
\hline
       & $1:10^0$     &  1.67e-07     &   359      \\
       & $1:10^1$     &  7.37e-07     &   442      \\
       & $1:10^2$     &  9.94e-07     &   723      \\
 SNNW-P& $1:10^3$     &  1.94e-06     &   809      \\
       & $1:10^4$     &  2.04e-06     &   773      \\
       & $1:10^5$     &  2.11e-06     &   766      \\
       & $1:10^6$     &  2.02e-06     &   765      \\
 \hline
       & $1:10^0$     &  1.65e-07     &   359      \\
       & $1:10^1$     &  7.14e-07     &   442      \\
       & $1:10^2$     &  1.03e-06     &   723      \\
 SNNW-G& $1:10^3$     &  2.05e-06     &   809      \\
       & $1:10^4$     &  2.07e-06     &   773      \\
       & $1:10^5$     &  2.11e-06     &   766      \\
       & $1:10^6$     &  2.07e-06     &   765      \\
 \hline
       & $1:10^0$     & 5.79e-07      &  2000      \\
       & $1:10^1$     & 1.96e-06      &  2000       \\
       & $1:10^2$     & 6.96e-06      &  2000       \\
 SNNW-R& $1:10^3$     & 1.28e-05      &  2000      \\
       & $1:10^4$     & 1.42e-05      &  2000       \\
       & $1:10^5$     & 1.52e-05      &  2000       \\
       & $1:10^6$     & 1.56e-05      &  2000       \\
\hline
\end{tabular}
\end{center}
\end{table}

\newpage
\section{Conclusion}
In this paper, we  present SNNW for solving the partial differential equation in weak form with high accuracy.
The basic idea of SNNW is to use some functions based on neural networks  as base functions to span a subspace, then find an approximate solution in this subspace.
Training base functions and finding an approximate solution can be separated, that is
different  methods can be  used to train these base functions, and different methods can also be used to find an approximate solution.
We train these base functions by minimizing the loss function, which can based on different form, such as the strong form and the weak form.
In this paper, we use PINN, DGM or DRM to train these base functions in the subspace, and use the Galerkin method  to solve the partial differential equation.
Our method can achieve high accuracy with low cost of training.
Numerical examples show that the cost of training these base functions
is low, and only several hundred to two thousand epochs are needed for most  tests.
The error of our method can  fall below the level of $10^{-7}$ for some tests.

\section*{Acknowledgements}
This work is partially supported by National Natural Science Foundation of China (12071045) and Funding of National Key Laboratory of Computational Physics.


\end{document}